\documentclass[12pt]{article}
\usepackage{amsthm,amsmath,amsfonts,amssymb}
\usepackage{graphicx}
\usepackage{natbib}
\usepackage{mathtools}
\usepackage{dsfont}
\usepackage{float}
\usepackage{xcolor}
\usepackage{subcaption}
\usepackage{hyperref}
\hypersetup{colorlinks=true, linkcolor=blue, urlcolor=blue, citecolor=blue}

\theoremstyle{plain}
\newtheorem{theo}{Theorem}
\newtheorem{lemm}{Lemma}

\newtheorem{prop}{Proposition}
\newtheorem{defi}{Definition}
\newtheorem{assu}{Assumption}

\theoremstyle{remark}

\newtheorem{rema}{Remark}

\renewcommand{\P}{{\mathbb P}}

\renewcommand{\hat}{\widehat}
\renewcommand{\tilde}{\widetilde}
\newcommand{\bm}{\boldsymbol}

\newcommand{\R}{{\mathbb R}}

\newcommand{\I}[1]{\mathds{1}_{\{#1\}}}

\newcommand{\grid}{\mathcal{G}}
\newcommand{\france}{\mathcal{F}}

\newcommand{\given}{\ |\ }
\newcommand{\norm}[1]{\left|\left| #1 \right|\right|}

\definecolor{R}{RGB}{255, 150, 0}
\definecolor{T}{RGB}{0, 100, 255}

\newcommand{\blind}{0}

\usepackage{xr} 
\begin{document}

\def\spacingset#1{\renewcommand{\baselinestretch}%
{#1}\small\normalsize} \spacingset{1}

\vspace{-15mm}
\if0\blind
{
\title{
    \bf
    Assessing the size of spatial extreme events using local coefficients based on excursion sets
    }
 
  \author{Ryan Cotsakis$^{1,\mathrm a,}$\thanks{
    The authors gratefully acknowledge that this work has been supported by the French government, through the 3IA C\^{o}te d'Azur Investments in the Future project managed by the National Research Agency (ANR) with the reference number ANR-19-P3IA-0002. French weather data based on the SAFRAN reanalysis model were made available by  Météo France within its cooperation framework with INRAE. }  \and Elena Di~Bernardino$^{2,\mathrm b}$\and Thomas Opitz$^{3,\mathrm c}$}
    \date{\small{ \vspace{0.2cm}
    $^1$ Expertise Center for Climate Extremes (ECCE), Faculty of Business and Economics (HEC) - Faculty of Geosciences and Environment, University of Lausanne, CH-1015 Lausanne, Switzerland,\\
$^\mathrm a$ ryan.cotsakis@unil.ch \\
$^2$Universit\'e C\^ote d’Azur, Laboratoire J.A. Dieudonn\'e,   UMR CNRS 7351, Nice, $^\mathrm b$Elena.Di\textunderscore Bernardino@univ-cotedazur.fr\\
    \vspace{0.2cm}
    $^3$Biostatistics and Spatial Processes, INRAE, Avignon, France; $^\mathrm c$Thomas.Opitz@inrae.fr\\
    \vspace{0.2cm}
    \today
    }}
  \maketitle
} \fi

\if1\blind
{
  \bigskip
  \bigskip
  \bigskip
  \begin{center}
    {\LARGE\bf 
    Assessing the size of spatial extreme events using local coefficients based on excursion sets}
    \end{center}
      \medskip
} \fi
\vspace{-5mm}
\begin{abstract}
 
Extreme events arising in georeferenced processes can take various forms, such as occurring in isolated patches or stretching contiguously over large  areas, and can further vary  with the spatial location and the extremeness of the events. We use excursion sets above threshold exceedances in data observed over a  two-dimensional grid of rectangular pixels to propose a general family of coefficients that assess spatial-extent properties relevant for risk assessment, and study five candidate coefficients from this family. These coefficients are defined locally and interpreted as a spatial distance from a reference site where the threshold is exceeded. We develop statistical inference and discuss robustness to boundary effects and resolution of the pixel grid. To statistically extrapolate coefficients towards very high threshold levels, we formulate a semiparametric model and estimate a parameter characterizing how coefficients scale with the quantile level of the threshold. The utility of the new coefficients is illustrated through simulated data, as well as in an application to gridded daily temperature in continental France. We find notable differences in estimated coefficient maps between climate model simulations and observation-based reanalysis.

\end{abstract}

\noindent%
{\it Keywords:} asymptotic independence; climate modeling; extreme value theory; spatial statistics; tail inference.
\vfill
 
\spacingset{1.9} 

\section{Introduction}
Assessing the spatial structure of extreme events is essential for understanding and mitigating complex environmental risks, such as widespread heatwaves, severe droughts, or heavy rainfall episodes. In modern environmental-data settings, data are often available  on fine spatial grids with several observations over time, \textit{e.g.}, from climate model output or remote sensing imagery. This increasingly common data format enables the extraction of global shape information of the regions exhibiting threshold exceedances—so-called excursion sets—and provides an opportunity to complement traditional pair-based extremal dependence measures typically applied to relatively small numbers of irregularly-spaced observation locations.  
The excursion-set-based approach offers several practical advantages. Threshold exceedances are binary indicators, making computations memory-efficient and robust to outliers or censoring of data above extreme levels.

Classical multivariate and spatial Extreme-Value Theory (EVT) focuses extensively on pairwise tail dependence measures and their asymptotic behaviors, \textit{e.g.}, the bivariate extremal coefficient \citep{Schlather2003} or the extremogram \citep{coles1999,Davis2009,Strokorb2015}. These classical measures quantify dependence strength between pairs of locations but cannot fully capture how entire regions may simultaneously exceed critical levels. Moreover, traditional parametric spatial  models, often relying on strong assumptions or complex latent-variable structures \citep[\textit{e.g.},][]{Huser2024}, can be computationally demanding and lack flexibility in capturing intricate spatial patterns. Recent literature suggests that  environmental data frequently exhibit decreasing extremal dependence at higher marginal quantiles, so that exceedances become more localized as the threshold grows \citep{tawn2018,huser2022,wadsworth2022}. This behavior is typical for  asymptotic independence. Formally, for a spatial stochastic process $X(\bm s)$ with marginal distributions $F_{\bm s}$, asymptotic independence  between two locations $\bm s_1,\bm s_2$ means that the conditional exceedance probability $\P(F_{\bm s_1}(X(\bm s_1))>p\mid F_{\bm s_2}(X(\bm s_2))>p)$ tends to zero as the probability threshold $p$ tends to one. Quantifying and inferring such localization patterns from finite data remains challenging.

In the setting where data are numerical and indexed by a two-dimensional rectangular pixel grid, we introduce a novel class of scalar coefficients that directly summarize the geometric properties of excursion set regions exceeding a high, possibly location-dependent threshold (typically, a marginal quantile for a  high probability level), conditioned on an exceedance at a given reference site. Each coefficient is constructed to have a natural interpretation as a spatial measure with units of distance. These coefficients distill complex spatial dependence structures into interpretable, local summaries of the spatial ``footprint'' of extreme events. Unlike methods assessing only pair-based coefficients, these shape-based coefficients intrinsically capture how extremes spread out (or contract) in space, and how this behavior changes as thresholds become more extreme.

The proposed \emph{local excursion-set coefficients} are designed to: (i) adapt naturally to nonstationarity and complex domain boundaries; (ii) remain relatively simple to estimate and interpret from finite samples of the discrete random excursion set; and (iii) reflect the asymptotic behavior of the underlying random field. By assessing how these coefficients change as the threshold increases, we can infer whether the spatial extent of extremes shrinks, remains stable, or even grows in size.

In the present work, rather than advocating for a single summary statistic, we study a family of   coefficients, each capturing different aspects of the geometry of exceedance regions. For instance, some coefficients measure how quickly one encounters a non-exceedance point as we move away from the conditioning location, while others gauge the size of the connected exceedance cluster or capture the shape complexity of excursion boundaries. We present five examples of such coefficients in detail, discuss their properties, and provide guidance on their interpretation. Through theoretical considerations and simulation studies, we show that these coefficients can be consistently estimated and can be used to test for asymptotic independence. 
Indeed, while these five coefficients differ in their practical behavior---some are more conservative in their risk assessment (\emph{i.e.}, leading to relatively larger coefficient values), others are more robust to boundary effects or discretization---they share a fundamental property: under broad conditions, they all evolve at the same rate as thresholds grow large. In other words, at extremely high levels, all five coefficients essentially convey equivalent information regarding the underlying tail dependence structure of the random field. This theoretical equivalence suggests that there is no unique, universally optimal coefficient; each can serve as a diagnostic tool, chosen to best fit the goals of the analysis and the nature of the considered data.

Compared to well-known existing measures of spatial tail dependence, our approach provides a different and complementary perspective. Focusing on global shapes can improve inference in cases where pairwise extremes are rare (\emph{e.g.}, due to a short observation period, low observation frequency, or a very high threshold) or where one aims to characterize how an entire region transitions from relatively large-scale to more localized extremes as thresholds increase. The resulting coefficients can be mapped across space to highlight local differences in the spatial organization of extremes. This could be particularly relevant for climate risk assessments and environmental applications, where understanding the spatial footprint of hazards under future, more severe scenarios is increasingly critical.

The remainder of the paper is organized as follows. Section~\ref{sec:coefficients} introduces the local excursion-set coefficients and their theoretical foundations. Section~\ref{sec:NumericalStudy} evaluates these coefficients through simulation studies. Section~\ref{sec:inference} discusses inference methods for their asymptotic behavior as thresholds increase. Section~\ref{sec:application} applies the proposed methodology to French temperature data. Finally, Section~\ref{sec:conclusion} provides a discussion of the results and outlines future research directions. Supplementary materials (provided separately) include an additional figure and table.

\section{Local excursion sets-based coefficients}
\label{sec:coefficients}

\subsection{Beyond pair-based dependence measures for extremes}
\label{listeNOEL}


In this section, we introduce a class of geometric coefficients that describe the spatial extent of excursion sets, conditioned on a reference location exceeding a threshold. These coefficients quantify different distance- and area-based aspects of extremal behavior. Their design is guided by the following criteria:

$\bullet$ \textit{Interpretability in terms of spatial distance.} Coefficients should relate directly to distances to capture the notion of \emph{extent}.

$\bullet$  \textit{Computability on large spatial grids.} Modern data sets often include thousands or even millions of grid points, which requires efficient and low-cost algorithms.

$\bullet$  \textit{Local definition in space.} Coefficients should be defined locally, allowing spatial heterogeneity (non-stationarity) in dependence structures to be captured and mapped. 

$\bullet$  \textit{Conservatism in risk assessment.} Underestimating the spatial extent of extremes can lead to poor risk management; coefficients should provide conservative summaries.

$\bullet$  \textit{Robustness to boundary effects.} The presence of domain edges should not strongly bias estimates.

$\bullet$  \textit{Robustness to discretization effects.} Coefficients should remain stable when the grid resolution changes, avoiding artifacts or bias introduced by discretization that is coarse relative to the scale of spatial dependence.

$\bullet$ \textit{Ability to capture asymptotic dependence or independence.} The spatial extent of co-occurring extremes may shrink (asymptotic independence) or remain extensive (asymptotic dependence) as thresholds grow.

Following these considerations, we next introduce in the next section a class of  coefficients (see Section \ref{notations}) and highlight their theoretical and practical properties (see Section  \ref{examples}).

\subsection{Notation and definitions}\label{notations}

 Let $(\Omega, \mathfrak{F}, \P)$ be a probability space, and let  
$X : \Omega \times \mathbb{R}^2 \to \mathbb{R}$ 
be a random field, indexed on a discrete grid. If data are missing at certain locations, the corresponding points are excluded from the grid. 

\begin{defi}[Discrete domain]
\label{gridG}
Let $\mathcal S \subset \mathbb{Z}^2$ be a finite subset containing the origin $\bm 0$. Define 
$\mathcal G = \{ M \bm s : \bm s \in \mathcal S \}$, 
where
\[
M =
\begin{pmatrix}
\Delta_x & 0 \\ 
0 & \Delta_y
\end{pmatrix}, \quad \Delta_x, \Delta_y > 0.
\]
The grid $\mathcal G$ represents regularly spaced observation points.
\end{defi}

In practice, \(\mathcal{G}\) represents raster cells (or image pixels) where \(\Delta_x \Delta_y\) is the area of a cell.  

Let $\|\cdot\|$ be a norm on $\mathbb{R}^2$, inducing a metric on $\mathcal G$. For a threshold function $u : \mathbb{R}^2 \to \mathbb{R}$, the binary indicator $X(\bm s) > u(\bm s)$ defines the excursion set as follows. 

\begin{defi}[Excursion set]\label{def:excursion_set}
The excursion set of $X$ above the threshold $u$ is
\[
\mathcal{E}^{\mathcal{G}, X}_u = \{ \bm s \in \mathcal{G} : X(\bm s) > u(\bm s) \}.
\]
\end{defi}

To summarize the size and shape of $\mathcal{E}^{\mathcal{G}, X}_u$, we define functions that map its conditional distribution---given $\bm s\in \mathcal{E}^{\mathcal{G}, X}_u$ for some $\bm s\in \grid$---to a nonnegative real value. Each of the functions satisfies the following distance-homogeneity property (\emph{i.e.}, rescaling the spatial domain by a positive factor will rescale the function value accordingly), so that their image may be interpreted as a distance, and thus prescribed any unit of distance. 

\begin{defi}[Distance-homogeneity]\label{def:distance-homogeneity}
Let $\mathcal H \subset \mathbb{R}^2$ be a finite set, $\bm s \in \mathcal H$, and let $\mu_{\mathcal H}$ be a probability measure on $\wp(\mathcal H)$, where $\wp$ denotes the power set.
A function $\theta: (\mathcal{H}, \bm s, \mu_{\mathcal H}) \mapsto x \in [0, \infty)$ satisfies distance-homogeneity if, for any $\lambda > 0$,
\begin{equation}\label{eqn:units_of_distance}
\theta(\lambda \mathcal{H}, \lambda \bm s, \mu_{\lambda \mathcal{H}}) = \lambda \theta(\mathcal{H}, \bm s, \mu_{\mathcal H}),
\end{equation}
where $\lambda \mathcal{H} = \{ \lambda \bm h : \bm h \in \mathcal{H} \}$ and $\mu_{\lambda \mathcal{H}}(\{\mathcal{A}\}) = \mu_{\mathcal{H}}(\{\lambda^{-1} \mathcal{A}\})$, $\forall \, \mathcal{A}\in \wp(\lambda\mathcal H)$.
\end{defi}

\begin{defi}[Local excursion-set coefficient]\label{def:class_theta}
For $\bm s\in \grid$ with $\P(\bm s \in \mathcal E_u^{\grid, X}) > 0$, a local excursion-set coefficient at $\bm s$ is defined by $\theta(\grid, \bm s, \mu)$, where $\theta$ is any function that satisfies Definition~\ref{def:distance-homogeneity}, and $\mu$ is the law of the conditional excursion set $\mathcal E_u^{\grid,X} \given \bm s\in \mathcal E_u^{\grid, X}$. To make the dependence on the threshold function $u$ and the law of the excursion set explicit in the notation, we will  write $\theta(\bm s; \mathcal{E}^{\mathcal{G}, X}_u) = \theta(\mathcal{G}, \bm s, \mu)$.
\end{defi}

The functions $\theta$ that satisfy Definition~\ref{def:class_theta} quantify the spatial extent of $\mathcal{E}^{\mathcal{G}, X}_u$ conditioned on $\bm s \in \mathcal{E}^{\mathcal{G}, X}_u$, as exemplified by the five coefficients $\theta_k$, $k = 1, \ldots, 5$, that we introduce in the following subsection.

\subsection{Examples and interpretations}\label{examples}
Before introducing the exemplary coefficients, connectivity (based on the adjacency of pixels in the \emph{rook} sense, \emph{i.e.}, along horizontal and vertical directions) and basic geometric quantities are introduced below. 

\begin{defi}[Connectivity]
For $\mathcal{H} \subseteq \mathcal{G}$, points $\bm a, \bm b \in \mathcal{H}$ are connected in $\mathcal{H}$, denoted $\bm a \stackrel{\mathcal{H}}{\sim} \bm b$, if there exists a sequence $(\bm s_1, \ldots, \bm s_n) \subset \mathcal{H}$ such that $\bm s_1 = \bm a$, $\bm s_n = \bm b$, and $\bm s_{j+1} - \bm s_j \in \{ (\pm \Delta_x, 0)', (0, \pm \Delta_y)' \}$ for $j = 1,\ldots,n-1$.
\end{defi}

\begin{defi}[Area and perimeter]
What we call the area of $\mathcal{H} \subseteq \mathcal{G}$ is given by 
$
A(\mathcal H) = |\mathcal H| \Delta_x \Delta_y.
$ The perimeter $P(\mathcal H)$ is
\[
P(\mathcal H) = \sum_{\bm s \in \mathcal H} \Big( 
\I{\bm s + (0, \Delta_y)' \in \mathcal{G} \setminus \mathcal{H}} 
+ \I{\bm s - (0, \Delta_y)' \in \mathcal{G} \setminus \mathcal{H}} 
\Big) \Delta_x 
+ \Big(
\I{\bm s + (\Delta_x, 0)' \in \mathcal{G} \setminus \mathcal{H}} 
+ \I{\bm s - (\Delta_x, 0)' \in \mathcal{G} \setminus \mathcal{H}} 
\Big) \Delta_y.
\]
\end{defi}

For $\bm s\in \grid$ and $r > 0$, let $\mathcal{B}(\bm s, r) = \{ \tilde{\bm s} \in \mathcal{G} : \| \tilde{\bm s} - \bm s \| \leq r \}$ denote the discrete ball of radius $r$ centered at $\bm s$ and truncated to $\mathcal{G}$.

The five coefficients (with possible values in $[0,\infty]$) are as follows (where any quantile function $F^{-1}$ is understood as the generalized inverse of $F$ in the usual way):

\begin{enumerate}
    \item 
    \textit{Lower extremal range quantiles.} Define $R_*(\bm s) = \inf \{ \| \tilde{\bm s} - \bm s \| : \tilde{\bm s} \in \mathcal{G} \setminus \mathcal{E}^{\mathcal{G}, X}_u \}$. Let $F_1(r;\bm s) = \P(R_*(\bm s) \leq r \given \bm s \in \mathcal{E}^{\grid,X}_u)$ be the   cumulative distribution function (CDF)  of $R_*(\bm s)$, conditional on $\bm s \in \mathcal{E}^{\mathcal{G}, X}_u$. Then, 
    \[
    \theta_1(\bm s; \mathcal{E}^{\mathcal{G}, X}_u, \alpha) = F_1^{-1}(\alpha; \bm s),\qquad \alpha\in (0,1).
    \]

    \item \textit{Upper extremal range quantiles.} Define $R^*(\bm s) = \sup \{ \| \tilde{\bm s} - \bm s \| : \tilde{\bm s} \stackrel{\mathcal{E}^{\mathcal{G}, X}_u}{\sim} \bm s \}$. Let $F_2(r;\bm s) = \P(R^*(\bm s) \leq r \given \bm s \in \mathcal{E}^{\grid,X}_u)$. Then,
    \[
    \theta_2(\bm s; \mathcal{E}^{\mathcal{G}, X}_u, \alpha) = F_2^{-1}(\alpha; \bm s),\qquad \alpha\in (0,1).
    \]
    \item \textit{Extremal confidence range.} For $r > 0$, let $F_3(r; \bm s) = \mathbb{E}[A(\mathcal{E}^{\mathcal{G}, X}_u \cap \mathcal{B}(\bm s, r)) \given \bm s \in \mathcal{E}^{\mathcal{G}, X}_u ] / A(\mathcal{B}(\bm s, r))$. Then,
    \[
    \theta_3(\bm s; \mathcal{E}^{\mathcal{G}, X}_u, \alpha) = \inf\{r > 0 : F_3(r; \bm s) < \alpha\},\qquad \alpha\in (0,1).
    \]
    \item \textit{Confidence region range.} For $\alpha\in (0,1)$, define
    \[
    \theta_4(\bm s; \mathcal{E}^{\mathcal{G}, X}_u, \alpha) = \sup_{\mathcal{D} \subseteq \grid} \{ \sqrt{A(\mathcal{D})/\pi} : \P(\mathcal{D} \subseteq \mathcal{E}^{\mathcal{G}, X}_u \given \bm s \in \mathcal{E}^{\mathcal{G}, X}_u) \geq 1 - \alpha \}.
    \]

    \item \textit{Area-perimeter ratio.} If $\P(P(\mathcal{E}^{\mathcal{G}, X}_u) > 0 \given \bm s \in \mathcal{E}^{\mathcal{G}, X}_u) > 0$, define
    \[
    \theta_5(\bm s; \mathcal{E}^{\mathcal{G}, X}_u) = \frac{8 \mathbb{E}[A(\mathcal{E}^{\mathcal{G}, X}_u) \given \bm s \in \mathcal{E}^{\mathcal{G}, X}_u]}{\pi \mathbb{E}[P(\mathcal{E}^{\mathcal{G}, X}_u) \given \bm s \in \mathcal{E}^{\mathcal{G}, X}_u]}.
    \]
\end{enumerate}

\begin{rema}\label{rem:immediate_interpretations}
The distribution function \(F_1(r; \bm s)\), which underpins \(\theta_1\), is studied in \cite{Ryan1} and shown to be closely related to the intrinsic volumes of excursion sets. The subset \(\mathcal{D}\) used in the construction of \(\theta_4\) is analyzed in \cite{bolin2015, bolin2018}. The coefficient \(\theta_5\) involves normalization by the factor \(\frac{8}{\pi} = \frac{4}{\pi} \times 2\), which serves two purposes: correcting for the discretization bias inherent in the perimeter length \(P\) \citep[see, \textit{e.g.},][]{cotsakis2022_1}, and aligning the area-perimeter ratio with the radius of an equivalent circular excursion set, specifically \(2 \times \frac{\pi r^2}{2\pi r} = r\). 
\end{rema}

\section{Evaluation of the considered coefficients}
\label{sec:NumericalStudy}
We now discuss how the coefficients \(\theta_1, \ldots, \theta_5\) provide complementary summaries of the spatial extent of extremes, balancing conservatism, computational efficiency, and interpretability. 

\textbf{Interpretability.} Following from Definition~\ref{def:distance-homogeneity}, each of the local excursion set coefficients is measured in units of distance. This means that each coefficient scales linearly (up to discretization error) with the amount of spatial rescaling.\\
The coefficients $\theta_1$ and $\theta_2$ are quantiles of mappings of realizations of $\mathcal E_u^{\grid, X}\given \bm s\in \mathcal E_u^{\grid, X}$ to real numbers that have clear geometric interpretations. Remark that $R_*(\bm s)$ is the supremum of all $r > 0$ that satisfy $\mathcal B(\bm s, r) \subseteq \mathcal E_u^{\grid, X}$, and $R^*(\bm s)$ is the infimum of all $r > 0$ that satisfy $B(\bm s, r) \supseteq \{\tilde {\bm s} \in \grid : \tilde{\bm s} \stackrel{\mathcal{E}^{\mathcal{G}, X}_u}{\sim} \bm s \}$ (see Figure~\ref{fig:upper_and_lower} in the Supplementary materials for a graphical illustration).\\
The extremal confidence range $\theta_3$ is the smallest radius $r$ around $\bm s$ such that one expects a proportion $\alpha$ of $\mathcal B(\bm s, r)$ to be contained in $\mathcal E_u^{\grid, X}$. We direct the reader to \cite{bolin2015} for an interpretation of the confidence region range $\theta_4$, which is the square root of the area of what they refer to as the \textit{excursion set}---taking on a different meaning from Definition~\ref{def:excursion_set}---corresponding to the conditional random field $X\given X(\bm s) > u(\bm s)$.  Finally, the area-perimeter ratio $\theta_5$ provides the radius of circular excursion regions when the grid spacings $\Delta_x$ and $\Delta_y$ are small in comparison to this radius, correcting for the overestimate of the perimeter of the circle as mentioned in Remark~\ref{rem:immediate_interpretations}.

\textbf{Locality.} Each of the coefficients should be interpreted as a summary of the spatial extent of extremes \textit{conditioned on a threshold exceedance at a site $\bm s$}, owing to Definition~\ref{def:class_theta}.

\textbf{Computability.} Each $\theta_k$ can be estimated via straightforward “plug-in” methods using observations, or Monte Carlo samples from a fitted or hypothesized random field model. Conditional expectations are estimated by sample means that include only realizations of $X$ satisfying the conditioning event. The quantities $R_*(\bm s)$ and $R^*(\bm s)$ can be identified by a breadth-first or depth-first search, starting from $\bm s$ in $\mathcal O(|\grid|)$ compute time.   The computation of $R^*(\bm s)$ for several $\bm s$ can be sped up by pre-computing the connected components of each realization of $\mathcal E_u^{\grid, X}$ and the vertices of the convex hulls of each connected component. Similarly, $R_*(\bm s)$ may be computed for several $\bm s$ by implementing the Fast Marching Method \citep{sethian1996} with the non-exceedances as seeds.

\subsection{Numerical study}\label{sec:NumericalStudy_setup}
\subsubsection{Setup}

We now illustrate the behavior of the five coefficients of Section~\ref{sec:coefficients} using simulated data from Gaussian random fields. Let $\mathcal{G}$ be the square lattice $\{-60,\dots,60\}^2$ with grid spacing $\Delta_x = \Delta_y =~1$.

\begin{figure}[H]
\centering
\includegraphics[width=0.42\linewidth]{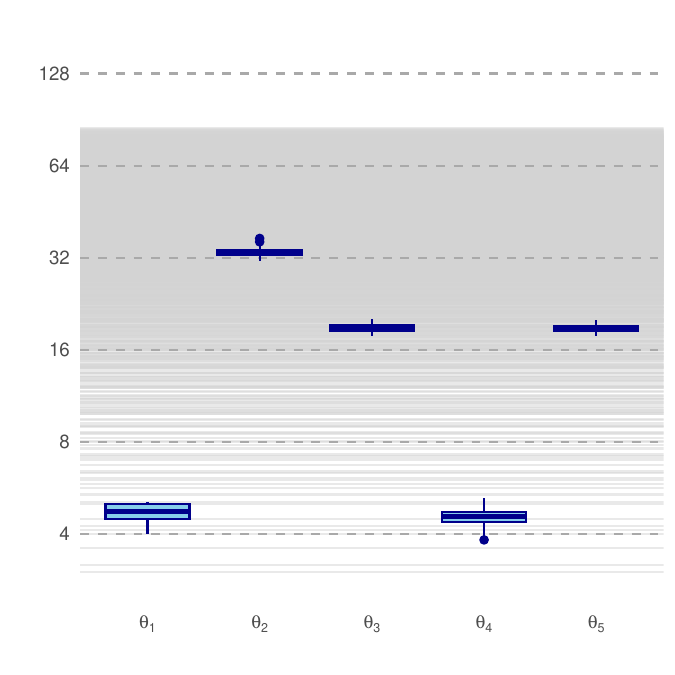}
\includegraphics[width=0.42\linewidth]{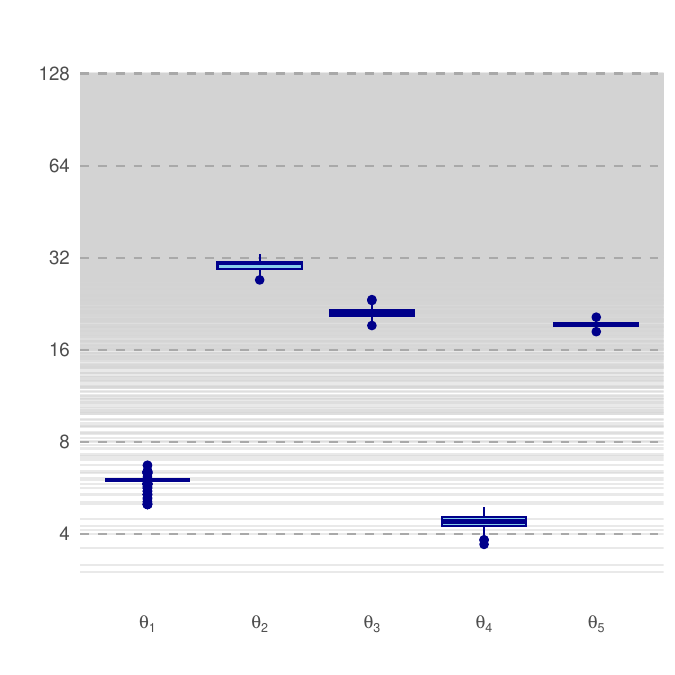}
\vspace{-0.7cm}\caption{Boxplots of $\theta_k(\bm s; \mathcal{E}^{\mathcal{G}, X}_u, 0.5)$, for $k=1, \ldots, 4$ and $\theta_5(\bm s; \mathcal{E}^{\mathcal{G}, X}_u)$ for $\bm s = \bm 0$ (left) and $\bm s = (60,0)'$ (right). A gray horizontal line is drawn  for each distance from $\bm s$ where there is a point in~$\grid$.}
\label{fig:boxplotsX}
\end{figure}

\begin{figure}[H]
\centering
\includegraphics[width=0.42\linewidth]{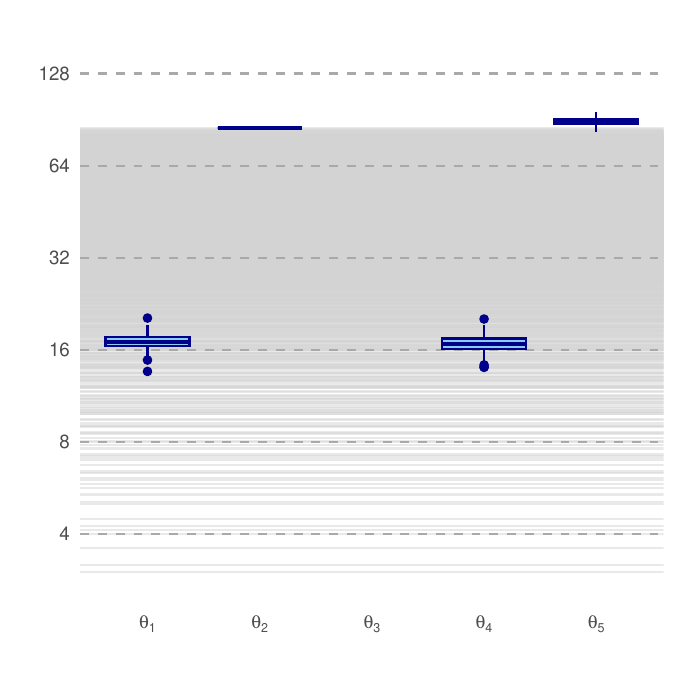}
\includegraphics[width=0.42\linewidth]{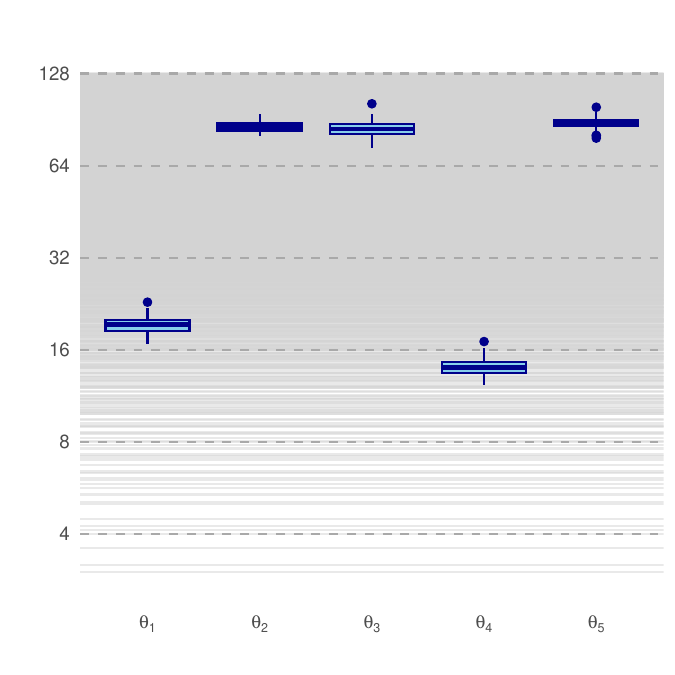}
\vspace{-0.7cm}\caption{Same as Figure~\ref{fig:boxplotsX} but for $\mathcal{E}^{\mathcal{G}, Y}_u$ instead of $\mathcal{E}^{\mathcal{G}, X}_u$.}
\label{fig:boxplotsY}
\end{figure}

We consider a zero-mean, unit-variance, stationary Gaussian random field $X$ on $\mathbb{R}^2$ with Matérn covariance: $\rho(\bm h) = \frac{2^{1-\nu}}{\Gamma(\nu)}\left(\frac{\sqrt{2\nu}\norm{\bm h}}{l}\right)^\nu K_\nu\left(\frac{\sqrt{2\nu}\norm{\bm h}}{l}\right),$ $\bm h\in \R^2$,
with regularity $\nu=2.5$ and range $l=30$, where $K_{\nu}$ is the modified Bessel function of the second kind. We also define $Y(\bm s) = X(\tfrac{\bm s}{4})$, which stretches the covariance range by a factor of 4.

We select $u$ to be the 99\% marginal Gaussian quantile, \emph{i.e.}\ $\mathbb{P}[X(\bm s) > u] = 0.01$. Here, since the considered model is stationary, $u$ is constant on $\mathcal{G}$. We simulate 500 i.i.d.\ realizations of $X$, conditionally on $X$ exceeding $u$ at the origin, $\bm s = \bm 0$. We obtain estimates of $\theta_k(\bm s; \mathcal E_u^{\grid,X}, 0.5)$, \textit{i.e.}, considering median levels with $\alpha=0.5$, for $k=1, \ldots, 4$, and of $\theta_5(\bm s; \mathcal E_u^{\grid,X})$, by approximating probabilities by sample frequencies and expectations by sample means. We repeat this process for $200$ bootstrap samples of the $500$ independent realizations of $\mathcal E_u^{\grid,X}\given \bm 0 \in \mathcal E_u^{\grid, X}$, and plot (on a logarithmic scale) the $200$ resulting estimates of each $\theta_k,\ k=1,\ldots,5$; results are reported in Figure~\ref{fig:boxplotsX} and interpreted in the following subsection. The analysis is repeated for the excursion set $\mathcal E_u^{\grid,Y}$ in place of $\mathcal E_u^{\grid,X}$ in Figure~\ref{fig:boxplotsY}. We also repeat the analysis for both fields using a boundary point $\bm s = (60, 0)'$ to test for biases induced by the edge of the domain. Table \ref{tab:boxplots} in the Supplementary materials summarizes the median and standard deviation of the logarithm of the estimates of  \(\theta_k\) across different ranges \(l\) in the Matérn covariance and locations \(\bm s\).  

\subsubsection{Results and discussion}

Several key observations emerge from the analysis of Figures~\ref{fig:boxplotsX}--\ref{fig:boxplotsY}.

\textbf{Conservatism.} Among the five coefficients, $\theta_2$ is the most conservative, consistently yielding the highest values. This makes it a relevant choice for risk assessment when understating the extent of an extreme event is problematic. However, this conservatism comes with trade-offs, particularly its sensitivity to boundary effects, as evidenced by the left plot of Figure~\ref{fig:boxplotsY}, where the median is consistently the maximum allowed value. In contrast, both $\theta_1$ and $\theta_4$ severely underrate the extent of an extreme event compared to the other coefficients and could lead to erroneous interpretation of the risk exposure. 

\textbf{Robustness to edge effects.} The coefficient $\theta_5$ stands out as the most robust to boundary effects, showing only negligible shifts when comparing results at the center of the domain $\bm{s} = \bm{0}$ and at the boundary $\bm{s} = (60,0)'$. In addition, this coefficient does not suffer from missing data at large distances from the conditioning location. By contrast, in Figure~\ref{fig:boxplotsY}, the left plot shows the sensitivity of $\theta_2$ and $\theta_3$ to domain truncation---$\theta_3$ is not even defined, as the estimated $F_3(\cdot;\bm 0)$ never drops below 0.5.

\textbf{Robustness to discretization.} As $Y$ has a correlation structure that is four times the extent of $X$, we expect that the quantities in Figure~\ref{fig:boxplotsY} should be four times those Figure~\ref{fig:boxplotsX}. Deviations from this factor of four are the result of spatial discretization error and statistical error---the latter can be accounted for by the width of the boxplots. The area-perimeter ratio $\theta_5$ performs well in this regard. Aside from edge truncation effects, $\theta_3$ also seems to reflect the factor of four  in the right-hand side of the two plots. In Figure~\ref{fig:boxplotsX}, on the left, we see that $\theta_1$ is comparable in size to $\Delta_x$ and $\Delta_y$, and its permissible values are quite limited due to the few number of points in $\grid$ with a corresponding magnitude of about 5.
\medskip

Overall, $\theta_5$ emerges as a well-balanced coefficient, performing reliably across all categories. It is moderately conservative, robust to edge effects, and scales consistently with changes in spatial dependence. While $\theta_2$ has an intuitive interpretation and is the most conservative of the coefficients, its sensitivity to edge effects limits its utility in quantifying asymptotic independence, especially when the domain size is comparable to the size of the extremes. For each of the coefficients, the low variability (on a log-scale) is indicative of the coefficient's ability to capture the extent of spatial extremes, and track the evolution as the threshold changes. We elaborate on this in the following section.

Code to reproduce Figures~\ref{fig:boxplotsX}--\ref{fig:boxplotsY} is provided at \url{https://github.com/pony-helping-treat/Local-Excursion-Set-Coefficients}.

If there is any specific limit or scaling behavior one wishes to highlight in real-world applications (\textit{e.g.}, reliance on large-domain asymptotics), care must be taken to assess possible discretization artifacts when $\Delta_x$ and $\Delta_y$ are not sufficiently small. This point will be further developed in Section~\ref{sec:inference} below.

\section{Inference for asymptotic behavior at increasingly extreme thresholds}
\label{sec:inference}

This section examines the asymptotics of the coefficients \(\theta_k\) as the threshold \(u\) grows. For each location $\bm s\in \grid$, let $F_{X(\bm s)}$ denote the CDF of $X(\bm s)$. Then, for each $p\in(0,1)$, define the quantile threshold function $u_p:\R^2 \to \R$ by
$u_p(\bm s)= F_{X(\bm s)}^{-1}(p)$.  

\subsection{Convergence of the considered coefficients}

 We  consider excursion sets where the dependence structure of the underlying random field is appropriately scaled by a function \(\sigma(p)\) so that a limit is obtained as $p\to 1$. In this way, the spatial extent of the extremes of random fields with no scaling must evolve like the reciprocal of $\sigma$. In what follows, we show that the coefficients $\theta_k$ can be used to uncover the asymptotic behavior of $\sigma(p)$ as $p\to 1$. Moreover, we propose a  model for $\sigma$ that gives rise to a more accurate inference strategy for the local excursion set coefficients $\theta_k$ at very extreme levels where excursions are scarce or yet unobserved in data.

\begin{defi}\label{Yfield} Let $\sigma: (0,1) \to (0,\infty)$ be a scaling function and let $p\in (0,1)$. Let us define
    \begin{equation}\label{eqn:setEp}\mathcal E^\grid_p = \left\{\bm s\in \grid : X\left(\frac{\bm s}{\sigma(p)}\right) > u_p\left(\frac{\bm s}{\sigma(p)}\right)\right\}.\end{equation} 
    The set $\mathcal E^\grid_p$ may be interpreted as the excursion set of the random field $Y:\Omega\times \R^2 \rightarrow \R$, defined by $Y(\bm s) = X\left(\frac{\bm s}{\sigma(p)}\right),\ \bm s\in \R^2$, at the threshold $u_p^Y:\R^2 \to \R$ defined by $u_p^Y(\bm s) = u_p\left(\frac{\bm s}{\sigma(p)}\right)$. For any fixed $p \in (0,1)$, the distribution of the conditional random set $\mathcal E_p^\grid\given \bm 0 \in \mathcal E_p^\grid$ on $\wp(\grid)$ is denoted $\nu_p$. 
\end{defi}

    Note that $\bm 0 \in \grid$ by construction in Definition~\ref{gridG}.

\begin{assu}\label{ass:rescalability} 
Suppose that $\sigma$, the scaling factor used to define $\mathcal E^\grid_p$ in~\eqref{eqn:setEp}, is chosen such that  $\forall \, \mathfrak{A}\in \wp(\wp(\grid))$  the probability measure $\nu_p$ on $\wp(\grid)$  satisfies
$$\nu_p(\mathfrak{A}) = \P\left(\mathcal E^\grid_p \in \mathfrak{A}\given \bm 0 \in \mathcal E^\grid_p\right) \xrightarrow[p \to 1]{} \nu(\mathfrak{A}),$$
with $\nu$  a non-degenerate limiting measure.
\end{assu}

\begin{rema}[On Assumption~\ref{ass:rescalability}] 
\label{rem:scaling_theory} 
One can adapt the proofs of Propositions 2 and 3 in \cite{Ryan1} in order to show that both Gaussian and regularly varying (RV) random fields satisfy Assumption~\ref{ass:rescalability}. Indeed, one can select the scaling function as $\sigma(p)= u_p(\bm 0)$ (\emph{i.e.}, linear behavior in the threshold) for stationary Gaussian random fields and $\sigma(p) = 1$ (\emph{i.e.}, constant behavior in the threshold) for RV random fields. For details on RV random fields, see \cite{Hult2005,Hult2006}. Coordinate rescalings of Gaussian processes that lead to nondegenerate extreme-value limits have also been studied in the context of max-stable processes \citep[Section~6,][]{Kabluchko2009}. 
\end{rema}

Each of the coefficients $\theta_k$, $k=1,\ldots 5$, is defined in terms of a location $\bm s\in\grid$ and the excursion set $\mathcal{E}^{\grid,X}_u$. By replacing $\mathcal{E}^{\grid,X}_u$ in the definition of our five coefficients with $\mathcal E^\grid_p$ in~\eqref{eqn:setEp}, 
we obtain (with some abuse of notation) the coefficients $\theta_k(\bm 0; \mathcal E^\grid_p, \alpha)$, for $k=1,\ldots,4,$ and $\theta_5(\bm 0; \mathcal E^\grid_p)$. Under Assumption \ref{ass:rescalability}, let $\mathcal E^\grid$ be the limiting  random element of $\wp(\grid)$ with distribution $\nu$, and define $\theta_k(\bm 0; \mathcal E^\grid, \alpha)$, for $k=1,\ldots,4,$ and $\theta_5(\bm 0; \mathcal E^\grid)$ accordingly.  

Finally, the following assumption ensures that each $\theta_k$,  $k = 1,\ldots,4$, is stable for small changes in the hyperparameter $\alpha$.

\begin{assu}\label{ass:quantile_convergence}
Fix $\alpha\in (0,1)$. Suppose that $\theta_k(\bm 0; \mathcal E^\grid, \tilde \alpha)$ is continuous in $\tilde \alpha$ in a small neighborhood of $\alpha$, for each $k=1,\ldots,4$.
\end{assu}

In the following, under Assumptions~\ref{ass:rescalability}--\ref{ass:quantile_convergence}, each of the five coefficients is studied for $p \to 1$.

\begin{prop}\label{prp:thetas_converge_too}
    Fix $\alpha\in (0,1)$. Let $\mathcal E^\grid_p$ be as in \eqref{eqn:setEp}.  Then, 
    \begin{itemize}
    \item under Assumption~\ref{ass:rescalability}, $\theta_5(\bm 0; \mathcal E^\grid_p) \xrightarrow[p\to 1]{}\theta_5(\bm 0; \mathcal E^\grid);$
    \item under Assumptions~\ref{ass:rescalability} and \ref{ass:quantile_convergence},
    $\theta_k(\bm 0; \mathcal E^\grid_p, \alpha) \xrightarrow[p\to 1]{}\theta_k(\bm 0; \mathcal E^\grid, \alpha)$, for $k = 1,\ldots,4$.
    \end{itemize} 
\end{prop}

\begin{proof}
    By the finiteness of $\grid$, the convergence of $\mathcal E^\grid_p  {\given \bm 0\in \mathcal E_p}$ to $\mathcal E^\grid$ implies the convergence of the conditional expectation of any real function of $\mathcal E^\grid_p$ to the expectation for $\mathcal E^\grid$. Moreover, the finiteness of $\grid$ implies that the set $\{\theta_k(\bm 0; \mathcal E^\grid_p, \alpha) : p \in (0,1)\}$, namely $\mathcal J$, has finite cardinality for all $k = 1,\ldots,4$. For $p\in (0,1)$, the coefficients $\theta_k(\bm 0; \mathcal E^\grid_p, \alpha)$ are monotonic in $\alpha$, and so we may partition $(0,1)$ into the intervals $\big(\{\alpha \in (0,1):\theta_k(\bm 0; \mathcal E^\grid_p, \alpha) = j\}\big)_{j\in \mathcal J}$. Under Assumption~\ref{ass:rescalability}, these intervals converge to $\big(\{\alpha \in (0,1):\theta_k(\bm 0; \mathcal E^\grid, \alpha) = j\}\big)_{j\in \mathcal J}$ as $p\to 1$, and under Assumption~\ref{ass:quantile_convergence}, $\alpha$ is an interior point of one of the limiting intervals.
\end{proof}

In the following, $\theta$ denotes any one of the coefficients $\theta_k$ for $k = 1,\ldots,5$, and the dependence on the hyperparameter $\alpha$ (in the case $k\neq 5$) is not written explicitly. The following lemma links the scaling behaviors of the coefficient, the excursion set and the grid.

\begin{lemm}\label{lem:rescaling_rules}
 The following equivalences hold for $p\in (0,1)$ under Assumption~\ref{ass:rescalability}:
\begin{enumerate}
       \item[(a)] $\lambda\,\theta(\bm 0; \mathcal E_p^{\grid}) = \theta(\bm 0; \lambda\,\mathcal E_{p}^{\grid})$, $\forall  \, \lambda > 0$. 
    \item[(b)] $\frac{1}{\sigma(p)}\,\mathcal E_p^\grid = \mathcal E_{u_p}^{\grid/\sigma(p), X}$ almost surely. 
    \item[(c)] $\theta(\bm 0;\mathcal E_p^\grid) = \sigma(p)\,\theta(\bm 0; \mathcal E_{u_p}^{\grid/\sigma(p), X})$.
\end{enumerate}
\end{lemm}
\begin{sloppypar}
\begin{proof} 
Item (a) is a direct application of Equation~\eqref{eqn:units_of_distance}, where $\theta(\bm 0; \lambda\,\mathcal E_{p}^{\grid})$ should be interpreted as the image of $(\lambda\grid, \bm 0, \mu)$ under $\theta$ where $\mu$ is the law of $\lambda\,\mathcal E_{p}^{\grid} \given \bm 0\in \mathcal E_{p}^{\grid}$ (see Definition~\ref{def:class_theta}). Item (b) holds almost surely  since
$\frac{1}{\sigma(p)}\mathcal E_p^\grid = \left\{\frac{\bm s}{\sigma(p)} : \bm s\in \grid, X\left(\frac{\bm s}{\sigma(p)}\right) > u_p\left(\frac{\bm s}{\sigma(p)}\right)\right\}
= \left\{\tilde {\bm s} \in \frac{1}{\sigma(p)}\grid: X(\tilde {\bm s}) > u_p(\tilde {\bm s})\right\} = \mathcal E_{u_p}^{\grid/\sigma(p), X}.
$ 
Item (c) follows directly from items (a) and (b) with $\lambda = 1/\sigma(p)$. 
\end{proof}
\end{sloppypar}

Results above show how, under Assumptions~\ref{ass:rescalability} and~\ref{ass:quantile_convergence},  one can choose an appropriate  scaling function $\sigma$ such that the limit coefficients    can be made  non-zero constants. 
 Alternatively,  one can choose the grid spacing parameters $\Delta_x$ and $\Delta_y$ sufficiently small, as we will show in the 
next section (see Assumption~\ref{ass:fine_grid} below).

\subsection{Inference of the scaling index}

Assumption~\ref{ass:sigma_regularly_varying} below is motivated by the behavior of Gaussian and regularly varying random fields and can be interpreted as a general model for the behavior of the appropriate scaling factor $\sigma(p)$ as a  function of the threshold $u_p$. 

\begin{assu}[Scaling index $\beta$]  \label{ass:sigma_regularly_varying}
    There exists $\beta\in \R$ such that for all $a > 0$, Assumption~\ref{ass:rescalability} holds for $\mathcal E_p$ in~\eqref{eqn:setEp} defined by a function $\sigma$ that satisfies
    \begin{equation}\label{eqn:sigma_regularly_varying}
        \lim_{x\to\infty}\frac{\sigma\left(1 - e^{-ax}\right)}{\sigma\left(1 - e^{-x}\right)} = a^\beta.
    \end{equation}
\end{assu}

\begin{rema}[On Assumption \ref{ass:sigma_regularly_varying}]Equation~\eqref{eqn:sigma_regularly_varying} is equivalent to asking that the  function $f:x\mapsto \sigma(1-e^{-x}),\ x\in (0,\infty)$ is regularly varying with index $\beta$. Since this corresponds to a marginal quantile at $p=1-e^{-x}$, we can $x$ as the marginal log-return time. The function $f$ is $\sigma$ composed with the CDF of a standard exponential random variable. Therefore, if $X(\bm 0)$ is standard exponential, then $u_p(\bm 0)$ corresponding to $p=1-e^{-x}$ is equal to $x$, $\forall\, x\in(0,\infty)$. Thus, in this case, $f$ maps $u_p(\bm 0)$ to $\sigma(p)$, relating the threshold to the appropriate scaling factor. Analogous to considerations in Remark \ref{rem:scaling_theory}, notice that Assumption \ref{ass:sigma_regularly_varying} is satisfied by  Gaussian random fields  with $\beta = 1/2$ and RV random  fields with $\beta = 0$. 
\end{rema} 

\begin{assu}\label{ass:fine_grid}
    Let $p \in (0,1)$. For any $\lambda, \epsilon > 0$, suppose that the grid $\grid$ may be chosen such that
    \begin{enumerate}
        \item[(a)]
        $\left|\theta(\bm 0;  \mathcal E_{u_{p}}^{\lambda\grid, X}) - \theta(\bm 0;  \mathcal E_{u_{p}}^{\grid, X})\right| < \epsilon$,
    \end{enumerate}
    and that under Assumption~\ref{ass:rescalability},  the scaling function $\sigma$ may be  {normalized}
    such that
    \begin{enumerate}
        \item[(b)] $\theta(\bm 0; \mathcal E^\grid) > 0$.
    \end{enumerate}
\end{assu}

\begin{theo}\label{thm:regression}
Suppose that for all $p_1,p_2\in (0,1)$ with $p_1 < p_2$, there exists a corresponding $\grid$ satisfying    Assumptions~\ref{ass:rescalability}--\ref{ass:fine_grid}. 
Then, for any $\epsilon >0$ and $a> 1$, there exists $p_0 \in (0,1)$  such that for all $p_1,p_2 \in (p_0,1)$ with $\log(1 - p_2) = a\log(1 - p_1)$, it is possible to choose $\grid$ such that
\begin{equation*}
        \left|-\beta - \frac{\log\theta(\bm 0; {\mathcal E}_{u_{p_2}}^{\grid, X}) - \log\theta(\bm 0; {\mathcal E}_{u_{p_1}}^{\grid, X})}{\log x_2 - \log x_1} \right| < \epsilon, 
    \end{equation*}   
where $x_i = -\log(1-p_i)$ for $i=1,2$, and $\beta$ is the scaling index in~\eqref{eqn:sigma_regularly_varying}.
\end{theo}

\begin{proof}[Proof of Theorem~\ref{thm:regression}]
Assumption~\ref{ass:fine_grid}~(b) together with Item~(c) in Lemma~\ref{lem:rescaling_rules} and Proposition~\ref{prp:thetas_converge_too} imply
$\log\sigma(p) + \log \theta(\bm 0;  {\mathcal E}_{u_p}^{\grid/\sigma(p), X}) \xrightarrow[p\to 1]{} \log \theta(\bm 0; \mathcal E^\grid).$ Therefore,
\begin{equation}\label{eqn:differences_logs_to_0}
    \left|\left(\log\sigma(p_1) + \log \theta(\bm 0; {\mathcal E}_{u_{p_1}}^{\grid/\sigma(p_1), X})\right) - \left(\log\sigma(p_2) + \log \theta(\bm 0; {\mathcal E}_{u_{p_2}}^{\grid/\sigma(p_2), X})\right)\right|
\end{equation}
can be made arbitrarily close to 0 for sufficiently large $p_0$.
Then, by adding and subtracting $\log\theta(\bm 0; {\mathcal E}_{u_{p_1}}^{\grid/\sigma(p_2), X})$ to~\eqref{eqn:differences_logs_to_0}, we have
\begin{align}\label{eqn:both_growthand_discretization}
    \log\sigma(p_1) - \log\sigma(p_2) \approx & \log\theta(\bm 0; {\mathcal E}_{u_{p_2}}^{\grid/\sigma(p_2), X}) - \log\theta(\bm 0;  {\mathcal E}_{u_{p_1}}^{\grid/\sigma(p_2), X})\\
    &+ \left(\log\theta(\bm 0; {\mathcal E}_{u_{p_1}}^{\grid/\sigma(p_2), X}) - \log\theta(\bm 0; {\mathcal E}_{u_{p_1}}^{\grid/\sigma(p_1), X})\right),\nonumber
\end{align}
where $\approx$ means that the two expressions can be made arbitrarily close to the same value by choice of $p_0$. Moreover, by Assumption~\ref{ass:sigma_regularly_varying}, we have
\begin{equation*}
    \log\sigma(p_1) - \log\sigma(p_2) \approx -\beta\log a = -\beta(\log x_2 - \log x_1).
\end{equation*}
By Assumption~\ref{ass:fine_grid}~(a), $\grid$ may be chosen such that $\big|\log\theta(\bm 0; {\mathcal E}_{u_{p_1}}^{\grid/\sigma(p_2), X}) - \log\theta(\bm 0; {\mathcal E}_{u_{p_1}}^{\grid/\sigma(p_1), X})\big|$ is negligible in comparison to $\log x_2 - \log x_1$. Therefore, we may divide the right-hand side of~\eqref{eqn:both_growthand_discretization} by $\log x_2 - \log x_1$ to obtain
$$\frac{\log\theta(\bm 0; {\mathcal E}_{u_{p_2}}^{\grid/\sigma(p_2), X}) - \log\theta(\bm 0;  {\mathcal E}_{u_{p_1}}^{\grid/\sigma(p_2), X})}{\log x_2 - \log x_1}\approx -\beta.$$
Notice that, without loss of generality, we may assume that $\sigma(p_2) = 1$, since rescaling $\sigma$ and $\grid$ by the same positive constant maintains Assumptions~\ref{ass:rescalability}--\ref{ass:fine_grid}.
\end{proof}

Theorem~\ref{thm:regression} provides a method to infer the scaling index of the appropriate scaling factor $\sigma$. Indeed, by studying the evolution of any of the five considered coefficients as a function of the threshold $u$, we may deduce the scaling index $\beta$ as the slope of a linear regression on appropriately rescaled axes, which in effect, determines how the spatial extent of extreme events tends to decrease at high thresholds.

\begin{rema}[On grid spacing]
The result assumes that \(\Delta_x, \Delta_y\) are small enough to approximate continuous scaling behavior. Unsurprisingly,  if the grid is coarse, discretization effects may obscure the true asymptotics, leading to incorrect estimates of \(\beta\).
\end{rema}

\subsection{Estimation strategy}
\label{subsec:EstimationStrategy}

We present our general procedure to infer the coefficients at a location $\bm s$ by first empirically evaluating  them at several threshold levels \(u\) and then fitting a linear regression based on Assumption~\ref{ass:sigma_regularly_varying} and Theorem~\ref{thm:regression}. We illustrate the approach by generating 2000 independent realizations of a random field \(X\) (with $X$ and $\grid$ are as defined in Section~\ref{sec:NumericalStudy_setup}), each conditioned on the event \(\{X(\bm{0}) > 1.75\}\). For each realization, we compute the five coefficients
\[
\theta_k(\bm{0}; \mathcal{E}_u^{\mathcal{G}, X}, 0.5),\quad k=1,\ldots,4,
\quad\text{and}\quad
\theta_5(\bm{0}; \mathcal{E}_u^{\mathcal{G}, X}),
\]
evaluated at several threshold levels \(u\). Below, we detail each step of the procedure.

\textbf{i) Transforming thresholds and gathering estimates.}
We select a range of threshold levels \(u\) such that $\log(-\log(1-F_{X(\bm 0)}(u)))$ is spaced across a suitable range of values---sufficiently high to capture extreme behavior but not so high that estimates become unduly noisy or unstable. Let us denote $u_1,\ldots, u_J$ these thresholds. For each \(u\), we use all 2000 conditional realizations of \(X\) to compute $\hat{\theta}_k(\bm{0}; \mathcal{E}_u^{\mathcal{G}, X}, 0.5)$, $k=1,\ldots,4$, and $\hat{\theta}_5(\bm{0}; \mathcal{E}_u^{\mathcal{G}, X})$.
This yields a point estimate of each coefficient at each threshold $u_j$, for $j \in \{1, \ldots J\}$. 

\textbf{ii) Bootstrapping for variance.}
To quantify uncertainty in each \(\hat{\theta}_k\), we employ a nonparametric bootstrap:
\begin{enumerate}
  \item[a)] We draw a bootstrap sample of the 2000 conditional realizations with replacement, and compute all   \(\hat{\theta}_k\) at each $u_j$ on this bootstrap sample.
  \item[b)]  We repeat the above step 200 times to obtain a bootstrap distribution for each \(\hat{\theta}_k\), for each threshold $u_j$, for $j \in \{1, \ldots J\}$.
\end{enumerate}
From these bootstrap distributions, we estimate the variance of \(\log \hat{\theta}_k\) for each $k = 1,\ldots, 5$ and $u_1,\ldots, u_J$. These variances guide the weights in our subsequent linear regression and inform our confidence intervals.

\textbf{iii) Weighted linear regression in log scale.}
We examine how \(\hat{\theta}_k(\bm{0}; \mathcal{E}_u^{\mathcal{G}, X})\) changes with \(u\) by fitting a weighted linear regression. Concretely, using covariate-observation pairs
\[
\left\{ \left(\log(-\log(1-F_{X(\bm 0)}(u_j))),\;\log\hat{\theta}_k\big(\bm{0}; \mathcal{E}_{u_j}^{\mathcal{G}, X}\big)\right) \;:\; j=1,\dots,J \right\},
\]
we perform a weighted least-squares, linear fit, with the $j$th weight set to the inverse of the bootstrap variance of \(\log\hat{\theta}_k\big(\bm{0}; \mathcal{E}_{u_j}^{\mathcal{G}, X}\big)\). This downweights threshold levels whose coefficient estimates exhibit higher uncertainty.

\begin{figure}
    \centering
    \includegraphics[width=0.45\linewidth]{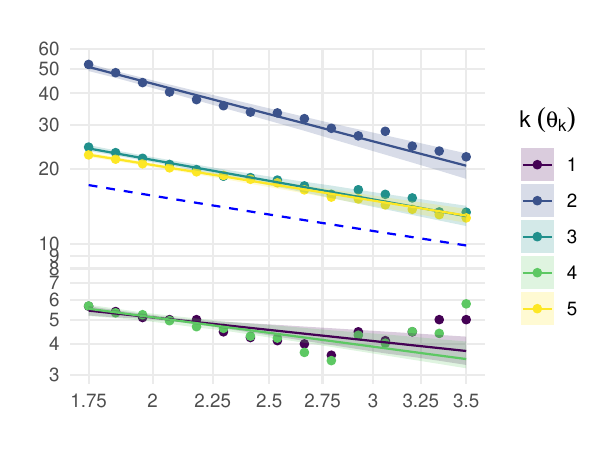}
    \hspace{.5cm}
    \includegraphics[width=0.45\linewidth]{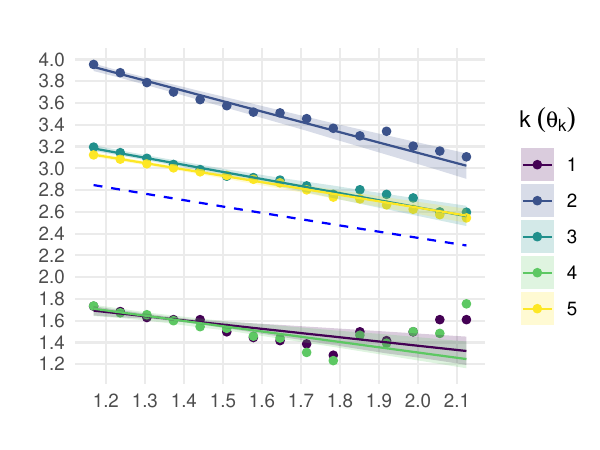}
    \put(-215,80){\rotatebox{90}{$\log\theta$}}
    \put(-440,90){\rotatebox{90}{$\theta$}}
    \put(-350,0){$u$}
    \put(-170,0){$\log(-\log(1-F_{X(\bm 0)}(u)))$}
    \caption{Illustration of the estimation strategy (see Section~\ref{subsec:EstimationStrategy}) carried out for the five considered local excursion set coefficients. The estimates $\hat\theta_k$ are displayed by dots, the fitted weighted linear regressions  by lines.  Twice the value of the expected area of $\{\bm s\in [-60,60]^2 : X(\bm s) > u\}$ divided by the expected length of $\{\bm s\in [-60,60]^2 : X(\bm s) = u\}$ is shown as a dashed blue curve. Both panels  display the same information but with different axes labels, so that the regression coefficients may be read off of the plot on the right.}
    \label{fig:simulation_gaussian}
\end{figure}

\textbf{iv) Bootstrapping for confidence intervals.}
The fitted regression line provides an “improved” estimate of \(\theta_k\) at \textit{any} threshold $u$. We then construct 95\% confidence intervals around this fitted line by repeating Step~3 for each of the 200 bootstrap samples of $\hat{\theta}_k$ computed in Step~2. In the resulting plots (see Figure~\ref{fig:simulation_gaussian}), we highlight these confidence intervals to give a visual assessment of the precision of our improved estimates. 

In Figure~\ref{fig:simulation_gaussian}, we observe that performing a linear regression on appropriately rescaled axes, in accordance with Theorem~\ref{thm:regression}, yields similar slope estimates across all coefficients. This consistency shows how each coefficient captures the same underlying degree of asymptotic independence. The regression lines decrease at a rate similar to the dashed blue curve that shows the theoretical value of $\theta_5$---calculated from the Gaussian kinematic formula \citep[][Theorem~15.9.5]{adler2007}---as the grid spacing approaches zero. Although the estimates of $\theta_5$ tend to slightly overestimate the theoretical values---attributed to the discrete nature of \(\mathcal{G}\)---the estimated slopes remain in close agreement with the theoretical expectations. Specifically, \(\theta_5\) exhibits a slope of \(-0.59\) with a 95\% confidence interval of \([-0.66, -0.52]\), which is near the theoretical slope that remains approximately \(-0.58\) in the range of the displayed axes (asymptotically, this slope tends to $-0.5$ as $u\to \infty$). Larger deviations are primarily observed for the smaller coefficients, where discretization effects play a more important~role.

\section{Application to French temperature data}\label{sec:application}

Temperature typically varies relatively smoothly in geographic space \citep{perkins2012}. Extreme events, such as heatwaves, can stretch over large areas. Assessing the spatial extent of such extremes, and how it changes with event frequency, is critical for understanding the compound risks they impose.

In this section, we illustrate and compare the proposed local excursion-set coefficients on two datasets of daily-average summer (June--August) temperatures over mainland France for their overlap period covering years 1970--2005. The majority of temperature extremes manifest in these summer months. Both datasets are defined on the same \(8\text{km}\times 8\text{km}\) grid, denoted \(\mathcal{F}\), in the Lambert-II projection. The two datasets under consideration are:

\begin{itemize}
    \item \textbf{Reanalysis}: Derived from the SAFRAN reanalysis model \citep{Vidal2010}, which integrates observational data with numerical weather prediction models to provide a comprehensive historical temperature record.
    \item \textbf{Simulation}: Generated for the CMIP5 historical experiment using the coupled IPSL-WRF global/regional climate model. This dataset is one of the reference models selected by the French weather service (Météo France) for studying climate change impacts (\url{http://www.drias-climat.fr/}).
\end{itemize}

We perform the same analyses for each of the two datasets, ensuring that any differences in the results can be attributed to either statistical uncertainties (which we assess) or more fundamental differences in the distributional properties of the two datasets.

For each location \(\bm{s} \in \mathcal{F}\), we compute the empirical \(p\)-quantiles \(u_p(\bm{s})\) for several \(p\) ranging from 0.86 to 0.99 using the $36$ summer seasons (approximately 3300 days) and we perform our analysis on the excursion sets 
$
\mathcal{E}_{u_p}^{\mathcal{F},X} = \left\{ \bm{s} \in \mathcal{F} \colon X(\bm{s}) > u_p(\bm{s}) \right\},
$
where \(X\) denotes either the Reanalysis or Simulation dataset. 

\begin{figure}[t]
    \centering
    \includegraphics[width=0.45\linewidth]{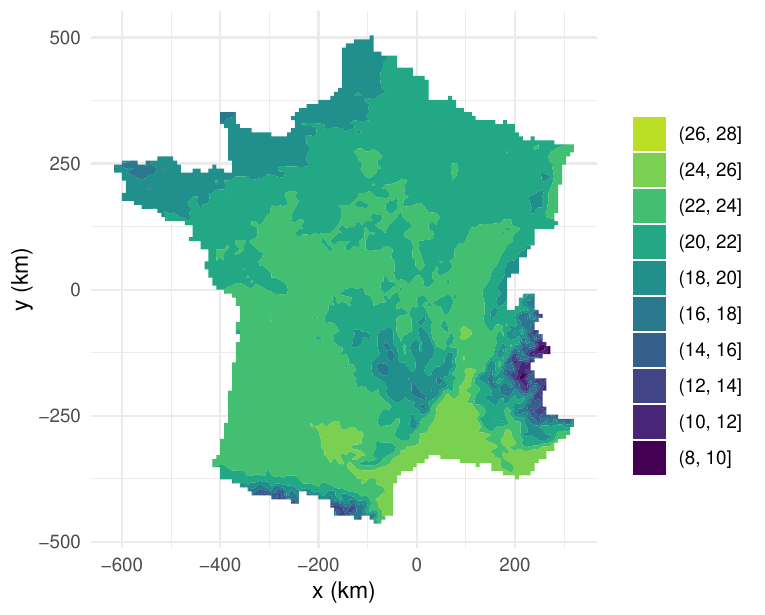}
    \put(-40,145){$T$ ($^\circ$C)}
    \includegraphics[width=0.45\linewidth]{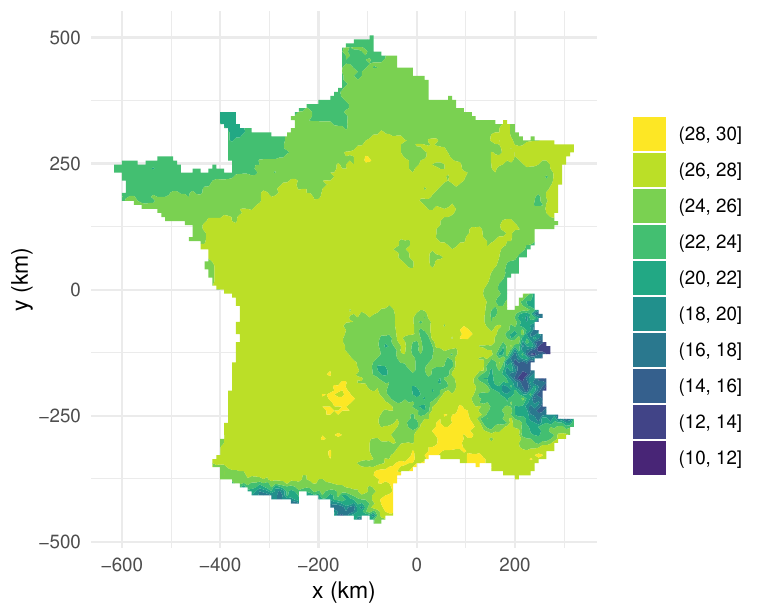}
    \put(-40,145){$T$ ($^\circ$C)}

\vspace{-0.25cm}    \caption{Estimated quantile thresholds \(u_p(\bm{s})\) for \(p = 0.86\) (left) and \(p = 0.99\) (right) across mainland France based on the Reanalysis dataset.}
    \label{fig:quantiles}
\end{figure}

Figure~\ref{fig:quantiles} illustrates the estimated 0.86 and 0.99 quantiles for the Reanalysis dataset, highlighting the spatial variability driven by topographic and climatic factors.

\begin{figure}[t]
\centering
\begin{subfigure}{.45\textwidth}
  \centering
  \includegraphics[width=0.9\linewidth]{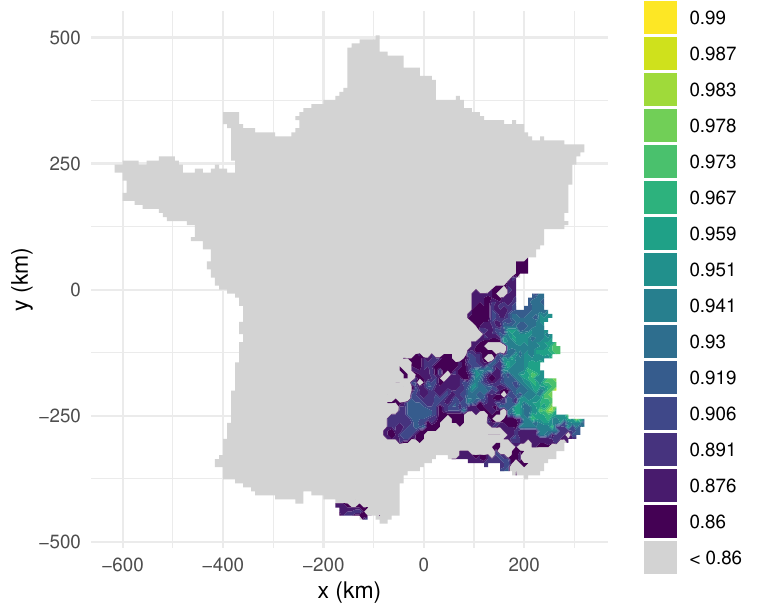}
  \put(-40,145){$p$}
  \caption{August 11, 1972}
  \label{fig:1972-08-11}
\end{subfigure}%
\begin{subfigure}{.45\textwidth}
  \centering
  \includegraphics[width=0.9\linewidth]{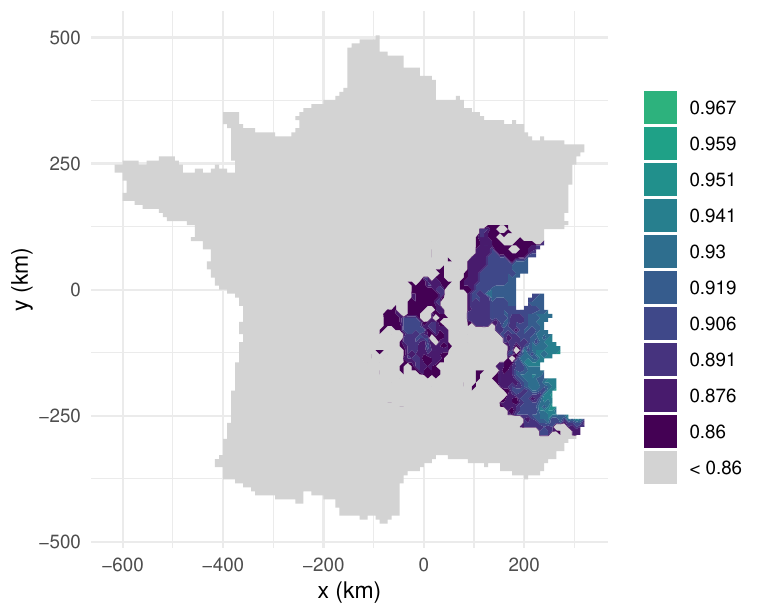}
  \put(-25,145){$p$}
  \caption{August 12, 1972}
  \label{fig:1972-08-12}
\end{subfigure}
\begin{subfigure}{.45\textwidth}
  \centering
  \includegraphics[width=0.9\linewidth]{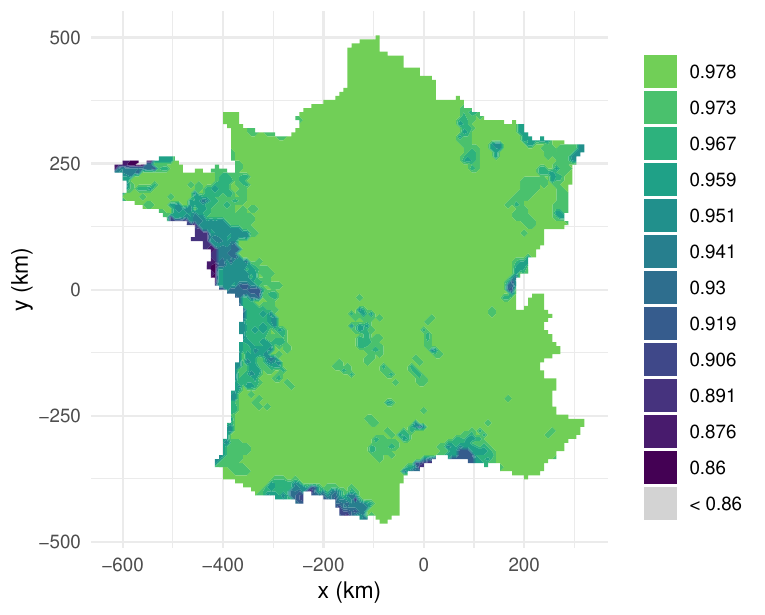}
  \put(-25,145){$p$}
  \caption{August 11, 2003}
  \label{fig:2003-08-11}
\end{subfigure}%
\begin{subfigure}{.45\textwidth}
  \centering
  \includegraphics[width=0.9\linewidth]{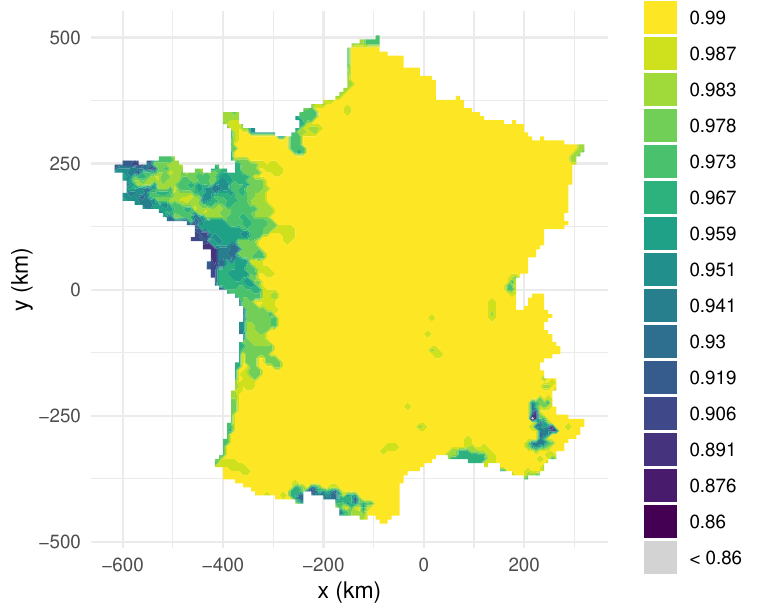}
  \put(-40,145){$p$}
  \caption{August 12, 2003}
  \label{fig:2003-08-12}
\end{subfigure}
\vspace{-0.25cm}\caption{$\mathcal E_{u_p}^{\mathcal F, X}$, for several $p$, observed on four different days. This corresponds to displaying the temperature over France on uniform margins, binned using the chosen values of $p$.}
\label{fig:temperatures}
\end{figure}

From the quantile maps for each dataset, we are able to compute the excursion set $\mathcal{E}_{u_p}^{\mathcal{F},X}$ for each quantile level $p$ on each day. We exhibit the excursion sets, with those for relatively larger $p$ nested within those for relatively smaller $p$, corresponding to the Reanalysis dataset on four different days in Figure~\ref{fig:temperatures}, with the two days in 2003 being part of one of the most severe heatwaves of recent decades. 

\subsection{Comparative analysis of local excursion set coefficients}

Following the estimation strategy outlined in Section~\ref{subsec:EstimationStrategy}, we compute the parameters of the linear models for the median upper extremal range \(\theta_2(\bm{s}; \mathcal{E}_{u_p}^{\mathcal{F},X}, 0.5)\) and the Area-Perimeter ratio \(\theta_5(\bm{s}; \mathcal{E}_{u_p}^{\mathcal{F},X})\) separately for each \(\bm{s} \in \mathcal{F}\). Given the temporal dependence inherent in daily temperature data, a stationary block bootstrap with a mean block size of 10 days is used when computing variances and confidence intervals \citep{Politis1994}.

\begin{figure}[t]
    \centering
    \includegraphics[width=0.6\linewidth]{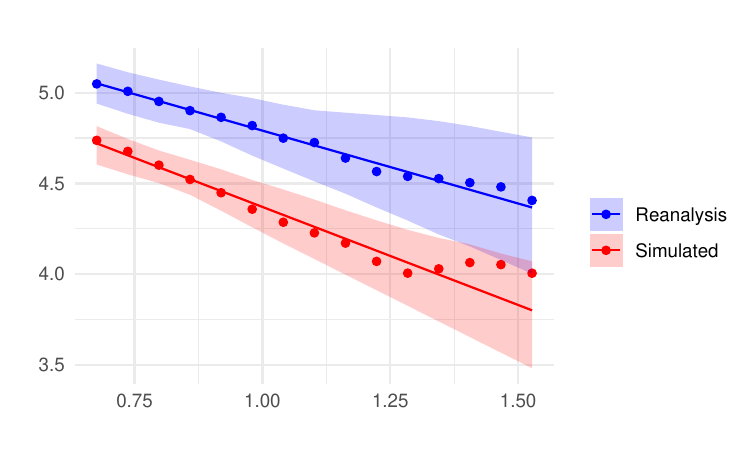}
    \put(-200,0){$\log(-\log(1-p))$}
    \put(-280,90){\rotatebox{90}{$\log\theta_5$}}
    \caption{Estimates of $\theta_5(\bm 0; \mathcal{E}_{u_p}^{\france, X})$ for both SAFRAN Reanalysis data (shown in blue) and for CMIP5 Simulation data (shown in red). A linear regression serves as an improved estimate of $\theta_5(\bm 0; \mathcal{E}_{u_p}^{\france, X})$, and the corresponding 95\% confidence interval (computed using stationary block bootstrap with mean block size of 10) is highlighted for each dataset.}
    \label{fig:single_location}
\end{figure}

The resulting fits of \(\theta_5(\bm{s}; \mathcal{E}_{u_p}^{\mathcal{F},X})\) for $\bm s = \bm 0$ (the origin of the coordinate system of Figure~\ref{fig:temperatures}) are shown in Figure~\ref{fig:single_location}. The difference between the model for the Reanalysis data and the Simulated data is statistically significant, as the confidence intervals (calculated from repeated linear fits using the bootstrap procedure) are non-overlapping for most of the considered region of threshold values. The two regression lines for $\log \theta_5$ are separated by at least $0.25$, such that estimated coefficients are at least $25\%$ higher for Reanalysis data. The possible discretization bias at high quantile levels $p$ (especially in Simulation data) is attenuated through the weighting procedure when calculating the linear fits. 
The discrepancies between the two datasets may stem from inherent differences in how the climate model (Simulation) represents physical processes governing extreme temperature events. Such biases could have profound implications for risk assessments and climate projections relying on model simulations. 

\subsubsection{Area-perimeter ratio \(\theta_5\)}
\begin{figure}
    \centering
    \includegraphics[width=0.48\linewidth]{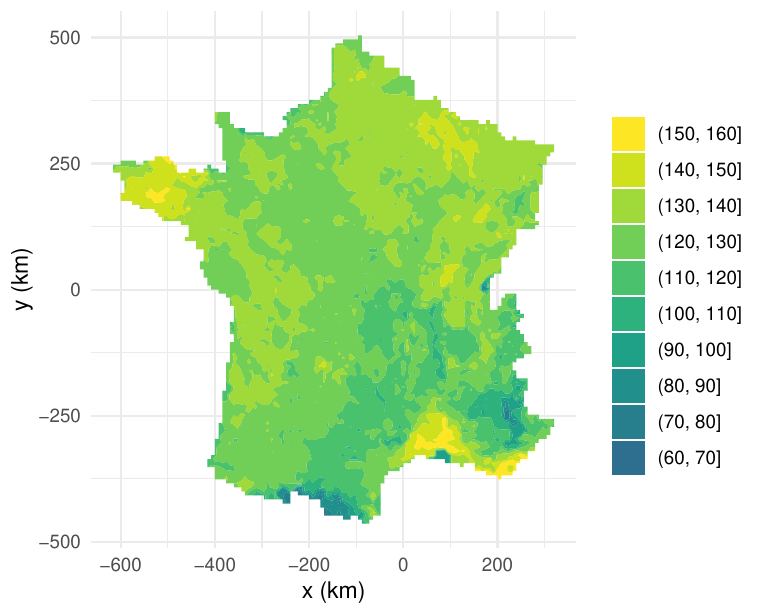}
    \put(-40,150){km}
    \includegraphics[width=0.48\linewidth]{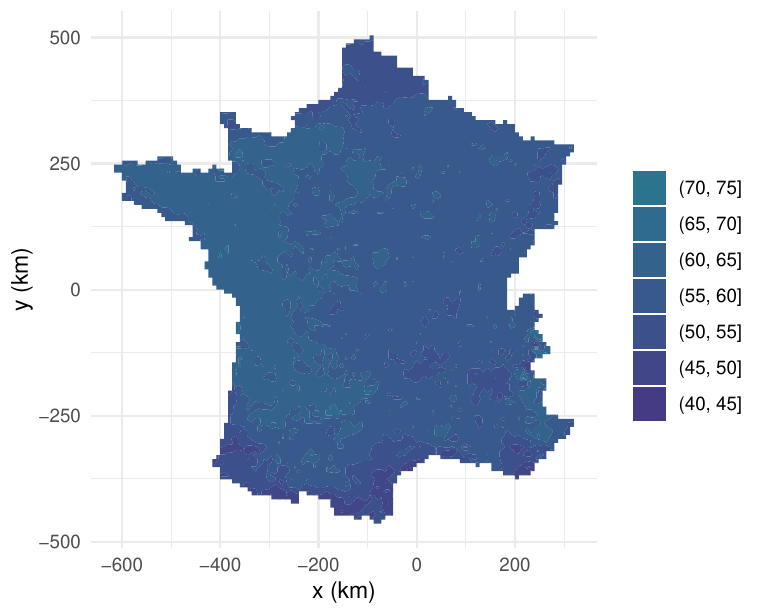}
    \put(-40,150){km}\\
    \includegraphics[width=0.48\linewidth]{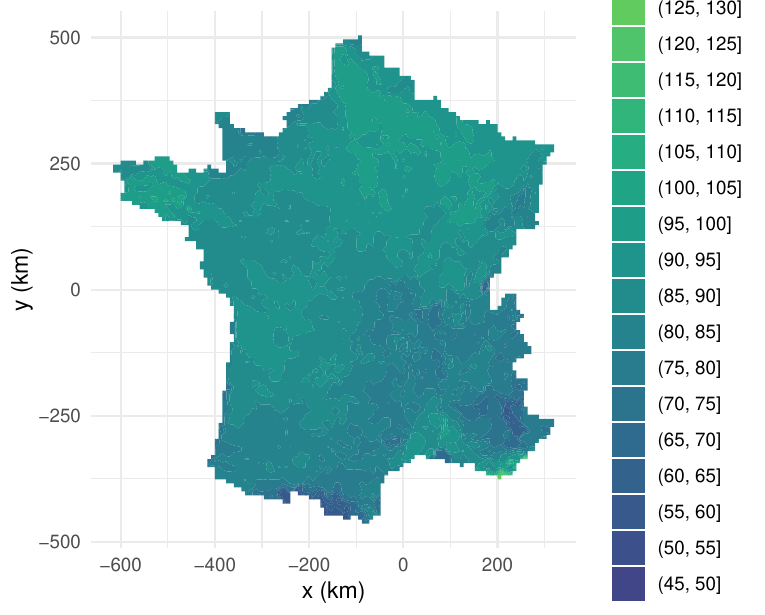}
    \put(-65,160){km}
    \includegraphics[width=0.48\linewidth]{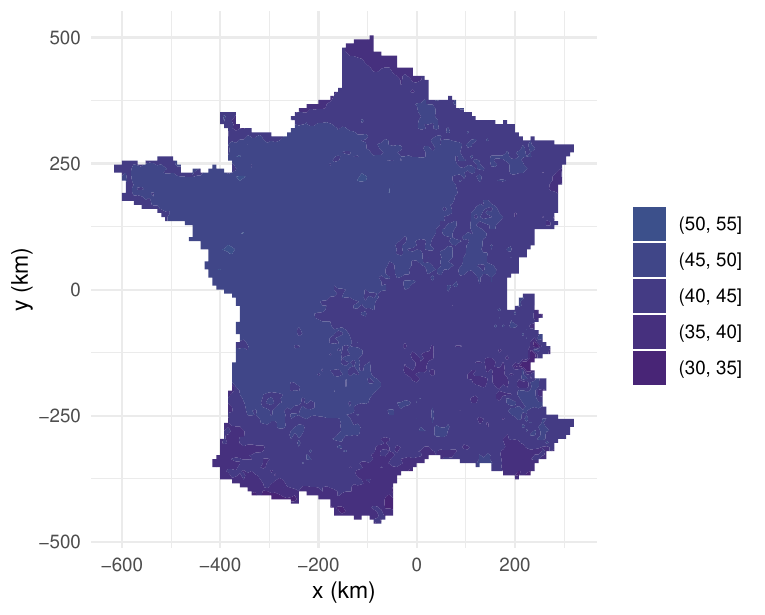}
    \put(-40,127){km}\\
    \includegraphics[width=0.48\linewidth]{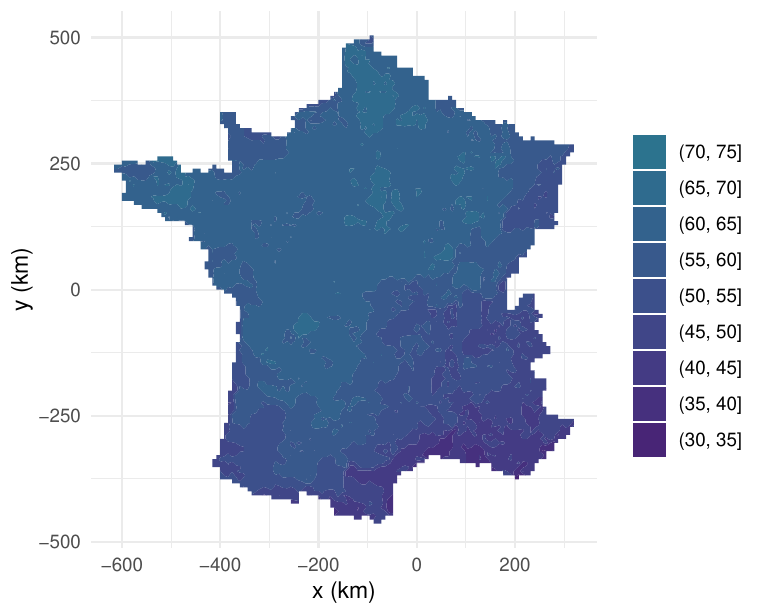}
    \put(-40,150){km}
    \includegraphics[width=0.48\linewidth]{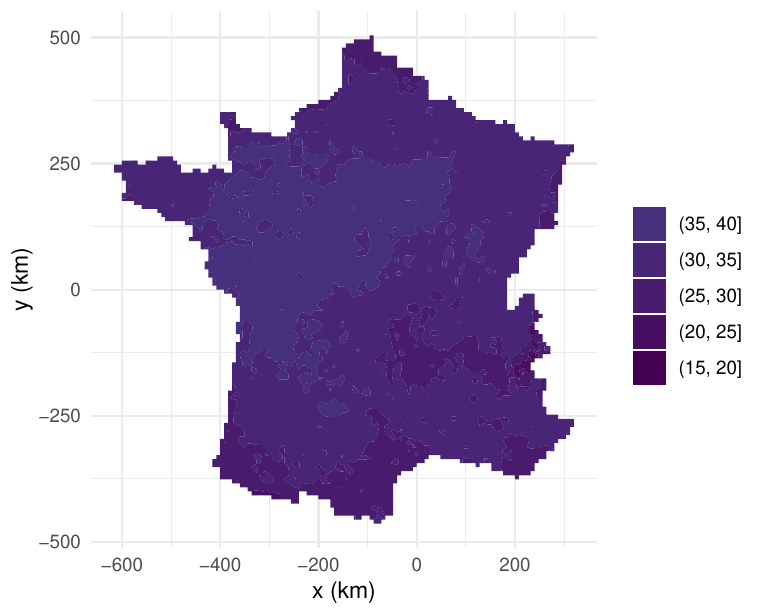}
    \put(-40,127){km}
    \caption{Estimated \(\theta_5(\bm{s}; \mathcal{E}_{u_{p}}^{\mathcal{F}, X})\) (measured in km) for the Reanalysis (left column) and Simulation (right column) datasets at \(p = 0.99\). The estimates are truncated at 160 km to enhance visualization. Top row: Upper end of the 95\% confidence interval. Middle row: The estimate $\hat\theta_5(\bm{s}; \mathcal{E}_{u_{p}}^{\mathcal{F}, X})$. Bottom row: Lower end of the 95\% confidence interval.}
    \label{fig:estimates_theta5}
\end{figure}

Figure~\ref{fig:estimates_theta5} presents a map of the estimated area-perimeter ratio \(\theta_5(\bm{s}; \mathcal{E}_{u_{p}}^{\mathcal{F}, X})\) at the 99th percentile threshold ($p = 0.99$) for both the Reanalysis and Simulation datasets. 
The Reanalysis dataset consistently exhibits larger \(\theta_5\) values compared to the Simulation dataset. Therefore, extreme temperature events in the Reanalysis data tend to cover larger spatial extents and form more coherent spatial clusters.
The differences in \(\theta_5\) between the two datasets are statistically significant across most grid points, as the estimated coefficients for the Reanalysis data fall outside the confidence interval for the Simulated data for all but $5$ southerly points in $\france$.
Both datasets show higher \(\theta_5\) values in regions with favorable conditions for heatwave persistence, such as low-lying plains.

\subsubsection{Upper extremal range median \(\theta_2\)}

Since \(\theta_2(\bm{s}; \mathcal{E}_{u_p}^{\mathcal{F},X}, 0.5)\) represents the median of \(R^*(\bm{s})\) conditioned on \(\bm{s} \in \mathcal{E}_{u_p}^{\mathcal{F},X}\), we leverage Generalized Additive Models (GAMs) to perform quantile regression (at the 50\% level) on \(\log R^*(\bm{s})\) whenever the logarithm is well-defined. Specifically, we use smooth splines to model spatial coordinates as covariates, capturing non-linear spatial trends, while imposing a linear dependence on \(\log(-\log(1-p))\) to reflect the relationship with the threshold level~\(p\). This modeling strategy would allow adding other covariates related to topography (\textit{e.g.}, elevation), or the time of observation to account for time trends. 

\begin{figure}
    \centering
    \includegraphics[width=0.48\linewidth]{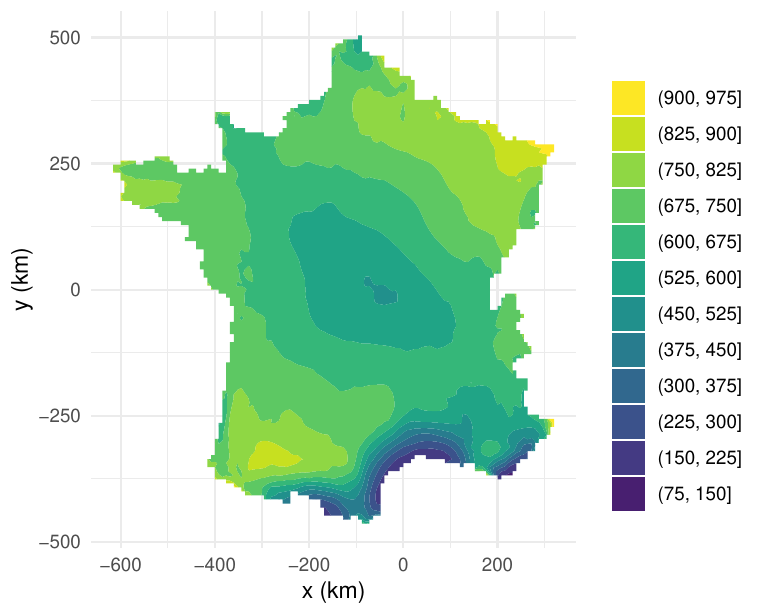}
    \put(-40,160){km}
    \includegraphics[width=0.48\linewidth]{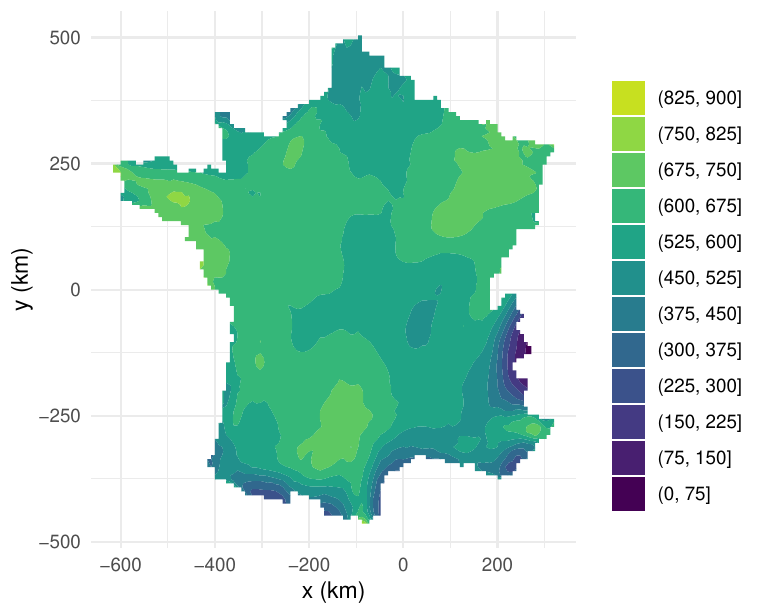}
    \put(-40,160){km}\\
    \includegraphics[width=0.48\linewidth]{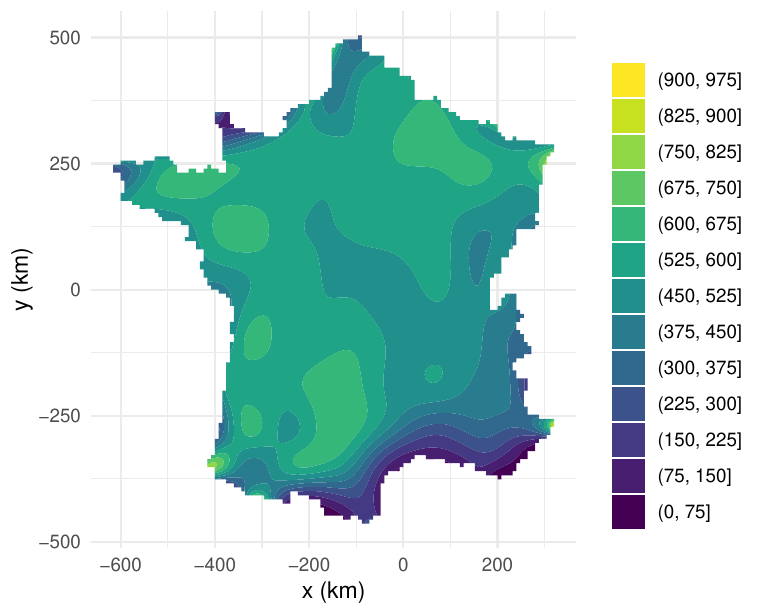}
    \put(-40,170){km}
    \includegraphics[width=0.48\linewidth]{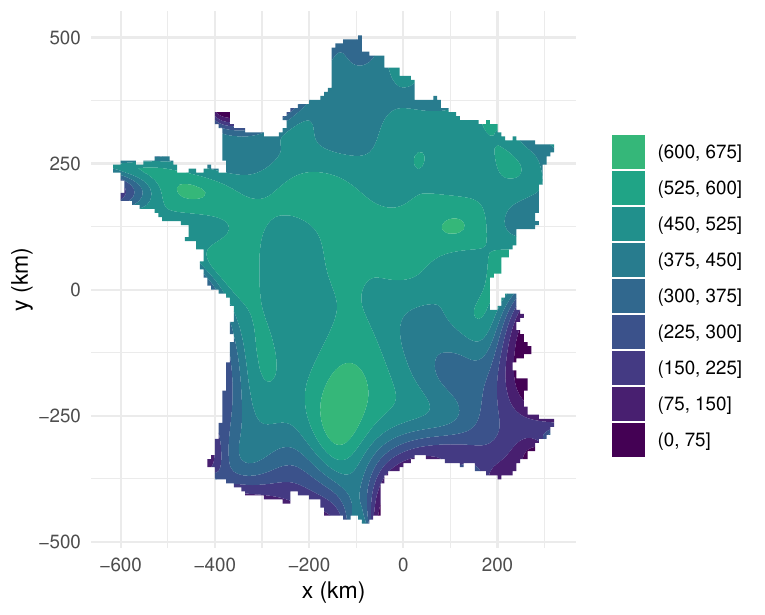}
    \put(-40,150){km}\\
    \includegraphics[width=0.48\linewidth]{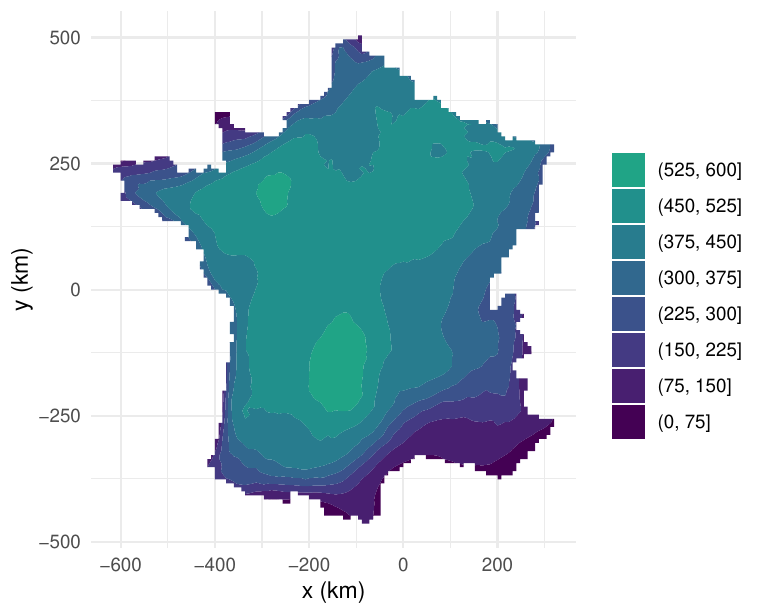}
    \put(-40,140){km}
    \includegraphics[width=0.48\linewidth]{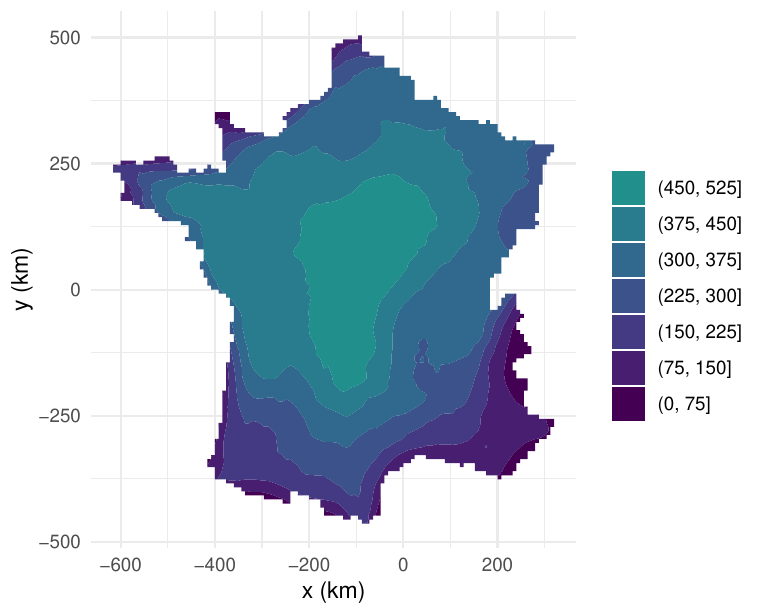}
    \put(-40,140){km}\\
    \caption{Estimated \(\theta_2(\bm{s}; \mathcal{E}_{u_p}^{\mathcal{F}, X}, 0.5)\) (measured in km) for the Reanalysis (left column) and Simulation (right column) datasets at \(p = 0.99\). Top row: Upper end of the 95\% confidence interval. Middle row: The estimate $\hat\theta_2(\bm{s}; \mathcal{E}_{u_p}^{\mathcal{F}, X}, 0.5)$. Bottom row: Lower end of the 95\% confidence interval. Left column: Reanalysis. Right column: Simulation.}
    \label{fig:estimates_theta2}
\end{figure}

Figure~\ref{fig:estimates_theta2} showcases the median upper extremal range \(\theta_2(\bm{s}; \mathcal{E}_{u_p}^{\mathcal{F},X}, 0.5)\) at the 99th percentile threshold of $X$. While both datasets demonstrate similar spatial patterns, the Simulation dataset generally displays smaller \(\theta_2\) values, reinforcing the observation of more localized extremes compared to the Reanalysis data.

The coefficient \(\theta_2(\bm{s}; \mathcal{E}_{u_p}^{\mathcal{F},X}, 0.5)\) provides a direct and interpretable measure of risk: it represents the median spatial extent of the connected component of the excursion set \(\mathcal{E}_{u_p}^{\mathcal{F},X}\) containing \(\bm{s}\), conditional on \(\bm{s} \in \mathcal{E}_{u_p}^{\mathcal{F},X}\). In this case study, larger \(\theta_2\) values for the Reanalysis dataset suggest that extreme temperature events are more widespread, while the smaller values for the Simulation dataset indicate a more localized spatial footprint of extremes.

However, due to the relatively large size of the extremes in this case study, \(\theta_2\) is subject to edge effects, particularly near the boundaries of \(\mathcal{F}\), where the connected component of the excursion set is truncated. The use of a GAM exacerbates these edge effects.  
 Consequently, while \(\theta_2\) is well-suited for assessing spatial risk and has an intuitive interpretation, it is not optimal for evaluating the degree of asymptotic independence, where boundary and discretization biases can distort the scaling behavior at high thresholds.

For judging the degree of asymptotic independence, \(\theta_5\) is more appropriate, as it is less prone to edge effects and exhibits less variability in its estimators at a given threshold (see Table~\ref{tab:boxplots}). Nonetheless, \(\theta_2\) remains a valuable complementary statistic, particularly for applications focused on risk quantification rather than asymptotic dependence analysis.

\subsection{Estimating the scaling index \(\beta\) using \(\theta_5\)}

To better understand how the spatial extent of extremes scales as thresholds become more extreme, we focus on the area-perimeter ratio \(\theta_5\) for estimating the scaling index \(\beta\). We follow the inference procedure  outlined in Section~\ref{subsec:EstimationStrategy}. The scaling index provides a quantitative measure of asymptotic independence, with lower \(\beta\) indicating weaker asymptotic independence (or stronger spatial dependence remaining at high thresholds).

\begin{figure}
    \centering
    \includegraphics[width=0.48\linewidth]{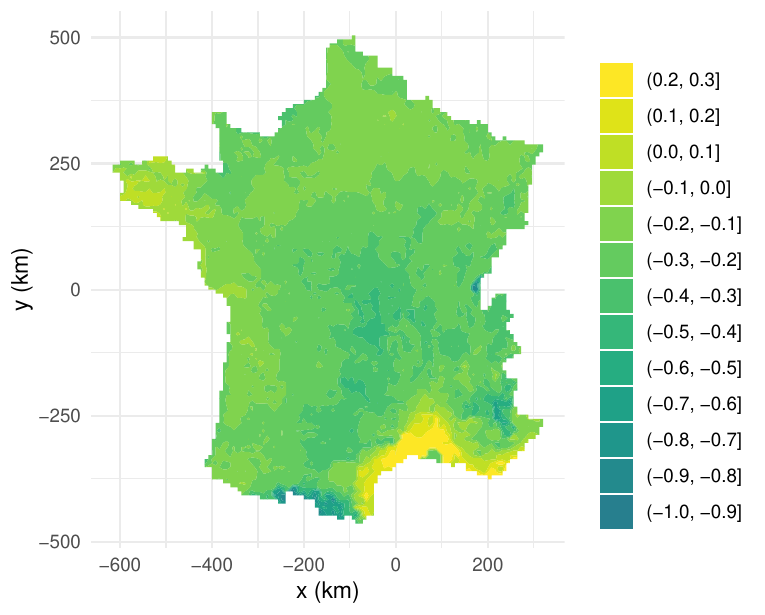}
    \includegraphics[width=0.48\linewidth]{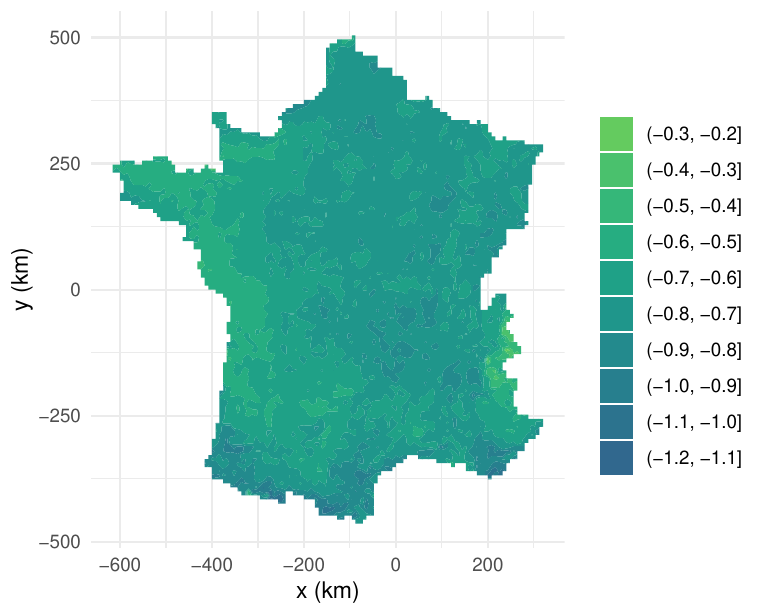}
    \includegraphics[width=0.48\linewidth]{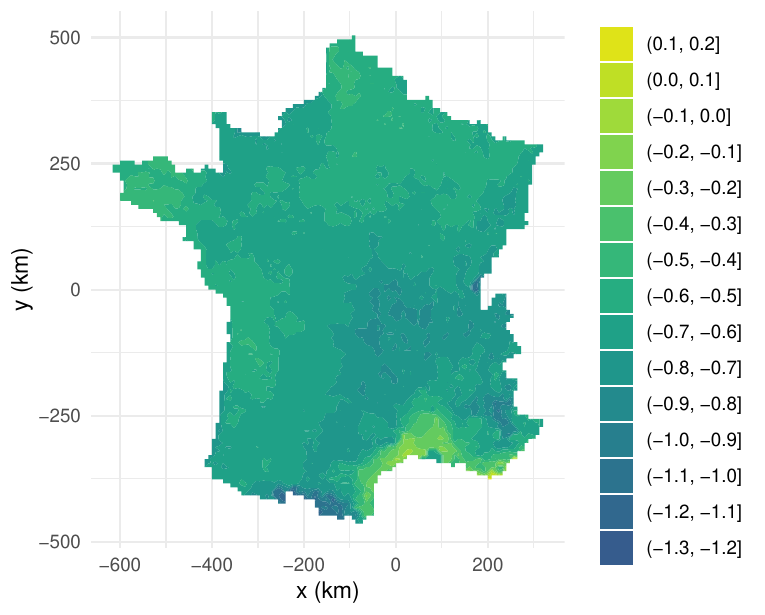}
    \includegraphics[width=0.48\linewidth]{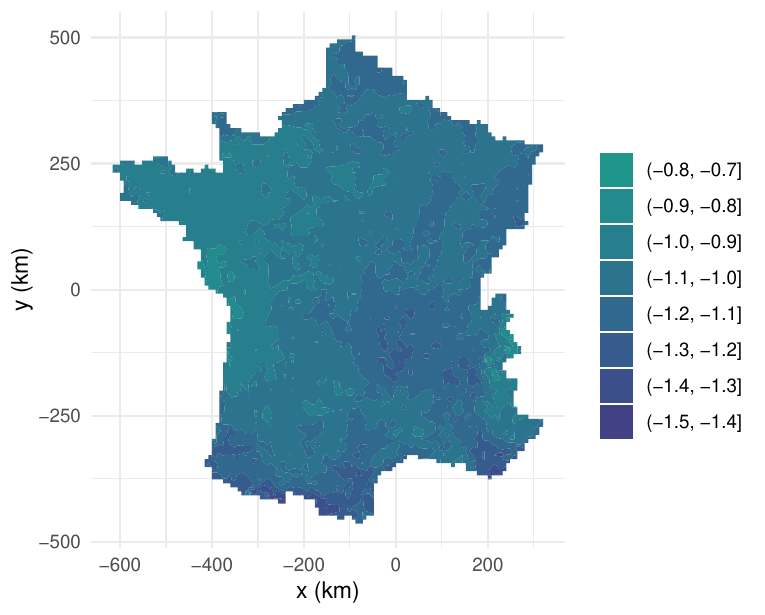}
    \includegraphics[width=0.48\linewidth]{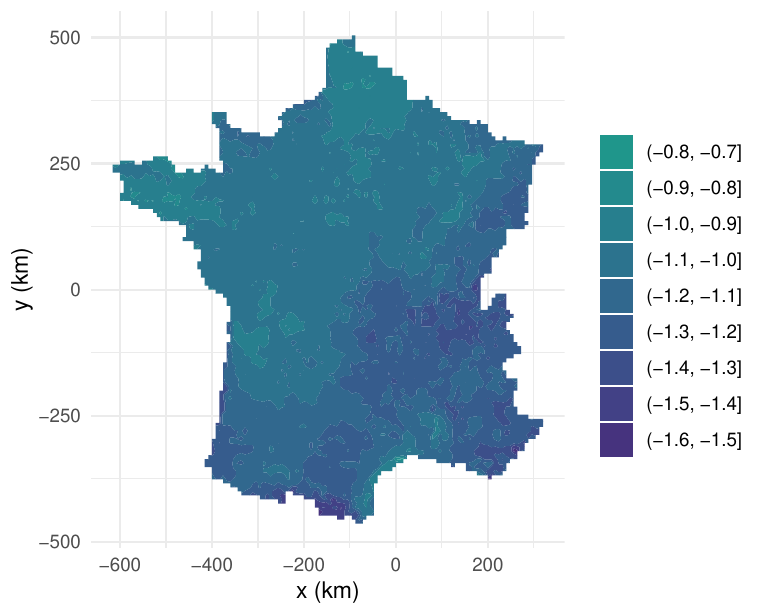}
    \includegraphics[width=0.48\linewidth]{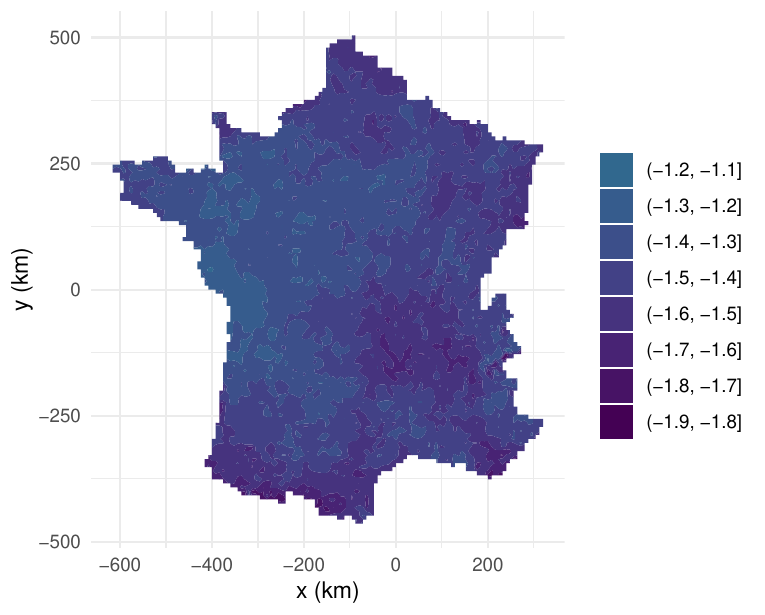}
    \caption{Spatial distribution of \(-\beta(\bm{s})\), the negative scaling index inferred from \(\theta_5(\bm{s}; \mathcal{E}_{u_p}^{\mathcal{F},X})\), across mainland France. Top row: Upper end of the 95\% confidence interval, indicating locations where asymptotic dependence cannot be rejected if the value is positive. Middle row: Estimated values of \(-\beta(\bm{s})\). Bottom row: Lower end of the 95\% confidence interval. Left column: Reanalysis dataset. Right column: Simulation dataset.}
    \label{fig:slopes}
\end{figure}

Figure~\ref{fig:slopes} maps the estimated \(\beta(\bm{s})\) for \(\bm{s} \in \mathcal{F}\). The Reanalysis data has \(\beta\) values of generally lower magnitude, especially in regions such as low-lying plains known to be prone to widespread heatwaves. In contrast, the Simulation data show higher \(\beta\) values uniformly, reflecting a tendency for extremes to become more localized as their intensity increases. 

We can test the null hypothesis of asymptotic dependence. It cannot be rejected at a given location \(\bm{s}\) if the confidence interval for \(-\beta(\bm{s})\) contains non-negative values, \textit{i.e.}, if its upper bound, as shown in the  top row of Figure~\ref{fig:slopes}, is positive. 
In Reanalysis data, asymptotic dependence cannot be rejected for areas near the Mediterranean coast in the South, where the confidence interval includes zero, suggesting weaker evidence of asymptotic independence in this region. However, for most other locations in mainland France, the confidence intervals exclude zero, allowing us to reject asymptotic dependence. 
In contrast, in Simulation data, asymptotic dependence is rejected everywhere.

\section{Discussion}
\label{sec:conclusion}

This work introduced local excursion-set coefficients to quantify the spatial extent and dependence of threshold exceedances in random fields. These coefficients provide interpretable measures of spatial dependence and enable new insights into the behavior of extremes, as demonstrated through theoretical analysis, simulations, and application to French temperature data. We have put focus on exceedances above relatively high marginal quantile levels and inference techniques for extrapolation towards very extreme levels, which is typically a challenging task. However, the coefficients we propose could also provide relevant information at intermediate levels, for example when the threshold corresponds to marginal medians where excursion sets would cover on average half of the study domain.
If data are not organized on a fine regular grid in metric units, they could be preprocessed using  interpolation methods (\textit{e.g.}, kriging, basis-functions techniques). 

The choice of an appropriate \(\theta_k\) depends on several factors. The scale of extremes relative to the spatial domain is crucial, ensuring minimal discretization and truncation error. Applications that prioritize risk aversion may favor more conservative coefficients like \(\theta_2\). Additionally, the resolution of the grid affects estimates, with coarser grids potentially introducing biases that influence the choice of \(\theta_k\). Finally, in scenarios where connectivity is critical, \(\theta_5\) has the drawback that it does not directly assess connected regions.
While this work focused on five illustrative examples, the general framework in Definition~\ref{def:class_theta} allows for a vast number of coefficients tailored to specific needs.

Strong links arise between the intrinsic volumes of unconditional excursion sets and the spatial extent of conditional extremes \citep{Ryan1}. This connection highlights the theoretical foundation of our coefficients and their potential for broader application in stochastic geometry. Our work illustrates  links between spatial extreme-value theory, stochastic geometry, and topological data analysis, and presents opportunities for further research to develop statistical tools for analyzing large-scale climate data.

The application to French temperatures revealed significant differences between Reanalysis and Simulation data. By estimating the scaling index \(\beta\) using \(\theta_5\), we showed that extremes in the Simulation dataset were more localized, indicating stronger asymptotic independence compared to the Reanalysis dataset. Additionally, we estimated \(\theta_2\), a very conservative measure of spatial risk, emphasizing its utility for practical risk assessment. 
Generalized Additive Models, as used with spatial position as covariate for \(\theta_1\) or \(\theta_2\)   could further incorporate other covariates like elevation, time, or large-scale atmospheric conditions \citep[similar to][]{Zhong2024,Koh2024} to infer more general nonstationary behavior. Moreover, applications to variables other than temperature could be explored. 

In future work, these coefficients could guide the construction of random field models exhibiting spatial nonstationarity in the dependence structure of extremes, such as determining the choice of noise distribution and of kernel bandwidths in process convolutions ({\it i.e.}, in fields constructed by smoothing pixel-based white noise), for example by using estimated maps of \(\theta_k\) as covariates.  They also have applications in generative models, including neural networks on grids, where they can capture spatial dependence in nonstationary generative processes like stable diffusion. 
Other useful methodological extensions would be to extend the definition of coefficients to spatiotemporal settings and multivariate processes, and develop connections to topological data analysis, ultimately enhancing the toolbox for studying climate extremes and informing decision-making.


{
\bibliographystyle{agsm} 
\bibliography{biblio}
\newpage

\section*{Supplementary materials}

\setcounter{section}{0}
\setcounter{figure}{9}
\setcounter{table}{0}
\setcounter{equation}{0}
\renewcommand{\thesection}{S\arabic{section}}
\renewcommand{\thefigure}{S\arabic{figure}}
\renewcommand{\thetable}{S\arabic{table}}

\def\spacingset#1{\renewcommand{\baselinestretch}%
{#1}\small\normalsize} \spacingset{1}
\spacingset{1.9}


To aid in the interpretation of our coefficients, Figure \ref{fig:upper_and_lower} illustrates the random variables \( R_{*}(\bm s) \) and \( R^{*}(\bm s) \), which are used in the construction of the coefficients \( \theta_1 \) and \( \theta_2 \), respectively. Specifically, \( R_{*}(\bm s) \) represents the distance from \(\bm s\) to the nearest non-exceedance point in \(\grid\), while \( R^{*}(\bm s) \) represents the distance from \(\bm s\) to the furthest connected exceedance point within the same connected component of the excursion set \(\mathcal{E}^{\mathcal{G}, X}_u\) as $\bm s$.

\begin{figure}[H]
    \centering
    \includegraphics[width=0.45\linewidth]{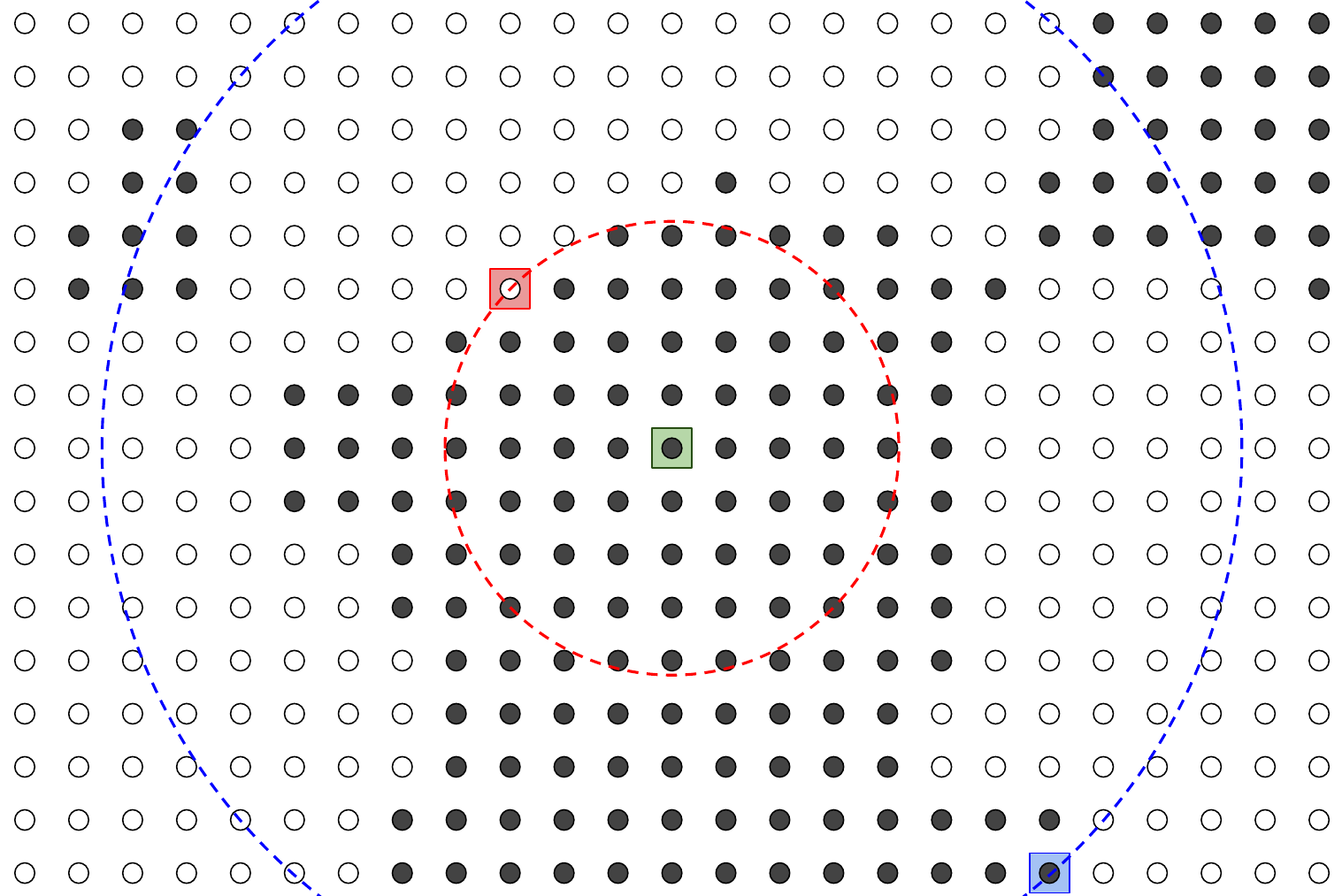}
    \caption{Illustration of excursion set geometry: the black points represent the excursion set \(\mathcal{E}^{\mathcal{G}, X}_u\) on a grid $\grid$. The green point is the chosen  reference location \(\bm s\). The red circle indicates \( R_{*}(\bm s) \)  and the blue one represents \( R^{*}(\bm s) \).}
    \label{fig:upper_and_lower}
\end{figure}


Table \ref{tab:boxplots} below summarizes the median and standard deviation of the logarithm of the estimated coefficients \(\theta_k\) across different ranges \(l\) and locations \(\bm s\).  
 
\begin{table}[H]
    \centering
    \small{\begin{tabular}{c||c|c||c|c}
         $\theta$ & $l$ (range) & $\bm s$ & Median of $\theta$ & S.d. of $\log\theta$ \\
         \hline
         $\theta_1$ & 30 & $\bm 0$ & 4.7 & 0.073 \\
          & 30 & $(60,0)'$ & 6.0 & 0.051 \\
          & 120 & $\bm 0$ & 17.0 & 0.057 \\
          & 120 & $(60,0)'$ & 19.3 & 0.063 \\
         \hline
         $\theta_2$ & 30 & $\bm 0$ & 33.3 & 0.029 \\
          & 30 & $(60,0)'$ & 30.7 & 0.035 \\
          & 120 & $\bm 0$ & 84.9 & 0.000 \\
          & 120 & $(60,0)'$ & 85.6 & 0.034 \\
         \hline
         $\theta_3$ & 30 & $\bm 0$ & 18.9 & 0.025 \\
          & 30 & $(60,0)'$ & 21.2 & 0.032 \\
          & 120 & $\bm 0$ & $\infty$ & NaN \\
          & 120 & $(60,0)'$ & 84.2 & 0.050 \\
         \hline
         $\theta_4$ & 30 & $\bm 0$ & 4.6 & 0.055 \\
          & 30 & $(60,0)'$ & 4.4 & 0.050 \\
          & 120 & $\bm 0$ & 16.7 & 0.060 \\
          & 120 & $(60,0)'$ & 14.0 & 0.061 \\
         \hline
         $\theta_5$ & 30 & $\bm 0$ & 18.8 & 0.021 \\
          & 30 & $(60,0)'$ & 19.3 & 0.021 \\
          & 120 & $\bm 0$ & 89.1 & 0.029 \\
          & 120 & $(60,0)'$ & 88.0 & 0.034
    \end{tabular}}
    \caption{Summary statistics of log-transformed coefficients   for different spatial ranges \(l\) and locations \(\bm s\). \(\bm s = \bm 0\) denotes the center of the grid, while \(\bm s = (60,0)'\) represents a boundary location. The range \(l\) indicates the spatial dependence scale in the simulation in Section~\ref{sec:NumericalStudy_setup}.}
    \label{tab:boxplots}
\end{table}
}
\end{document}